\documentclass{amsart}
\usepackage{amssymb, amsmath,bm}
\usepackage{color}
\usepackage{comment}
\usepackage{graphicx} 
\usepackage{multirow, array} 
\usepackage{float} 
\usepackage[english]{babel} 
\usepackage[normalem]{ulem}

\def\black{\color{black}}

\newtheorem{theorem}{Theorem}[section]

\newtheorem{remark}[theorem]{Remark}

\usepackage[dvipsnames]{xcolor}

\usepackage{pifont}

\newcommand{\Ricci}{\operatorname{Ric}}
\newcommand{\Riccicomp}[2]{\Ricci_{#2}^{\,\,#1}}
\newcommand{\RicciRiccicomp}[2]{\left(\Ricci^2\right)_{#2}^{\,\,#1}}

\def\span{\operatorname{span}}
\def\ve{\varepsilon}
\def\ip{\langle \cdot,\cdot \rangle}

\DeclareMathOperator{\ad}{ad}

\newcommand{\algsdR}{\mathbb{R}^3\rtimes \mathbb{R}}
\newcommand{\algsdH}{\mathfrak{h}_3\rtimes\mathbb{R}}

\newcommand{\alg}{\mathfrak{k}}

\newcommand{\Cotton}{C}

\usepackage{tikz,caption}
\usetikzlibrary{decorations.pathreplacing}

\begin{document}
\title[Four-dimensional pp-wave Lie groups and harmonic curvature]{Four-dimensional pp-wave Lie groups and harmonic curvature}\thanks{\emph{Data Availability Statement}: The data that support the findings of this study are available from the corresponding author upon reasonable request.}
\author[Garc\'ia-R\'io, Rodríguez-Gigirey, V\'azquez-Lorenzo]{E. Garc\'ia-R\'io, R. Rodríguez-Gigirey, R. V\'azquez-Lorenzo}
\address{EGR: Department of Mathematics, CITMAga, University of Santiago de Compostela, 15782 Santiago de Compostela, Spain}
\email{eduardo.garcia.rio@usc.es}
\address{RRG: IES Fuentesnuevas, 24411 Fuentesnuevas, Spain}
\email{rosalia.rodvil@educa.jcyl.es}
\address{RVL: IES de Ribadeo Dionisio Gamallo, 27700 Ribadeo,  Spain}
\email{ravazlor@edu.xunta.gal}
\thanks{\emph{Funding Statement}: Supported by projects  PID2022-138988NB-I00 (AEI/FEDER, Spain) and ED431C 2023/31  (Xunta de Galicia, Spain).}
	
\begin{abstract}
We determine all four-dimensional Lie groups which have harmonic curvature. As a consequence, a description of four-dimensional pp-wave Lie groups is obtained.
\end{abstract}
\maketitle

\section{Introduction}\label{se:introduction}

A pseudo-Riemannian manifold $(M,g)$ has \emph{harmonic Weyl curvature tensor} if its Weyl tensor is divergence-free, i.e., $\operatorname{div}W=0$, or equivalently the Cotton tensor vanishes. 
Since Einstein metrics have harmonic Weyl tensor,  the condition $\operatorname{div}W=0$ has been extensively investigated in order to extend geometric properties of Einstein metrics to more general contexts.
Trivially any pseudo-Riemannian manifold with parallel Ricci tensor has harmonic Weyl tensor, as well as locally conformally flat manifolds. It follows from \cite{datri, PS} that the converse is also true in the four-dimensional homogeneous Riemannian setting. Moreover, harmonicity of the Weyl tensor is equivalent to harmonicity of the curvature tensor, $\operatorname{div}R=0$, if the scalar curvature is constant.

Plane fronted waves with parallel propagation rays (pp-waves) are a special class of Lorentzian manifolds which provides the simplest mathematical family of exact solutions of Einstein's field equations. 
Plane waves are a particular class of pp-waves which includes gravitational plane waves as vacuum ones.
Since their introduction by Brinkmann in \cite{Bri}, pp-waves have appeared many times as solutions in problems motivated by either physical or purely geometrical questions. From the point of view of their curvature, pp-waves are the simplest non-flat Lorentzian manifolds since both the Ricci operator and the Weyl curvature operator are two-step nilpotent. Furthermore, while plane waves have harmonic Weyl tensor, pp-waves are not Cotton-flat in the generic situation.
 
Due to the special geometrical significance of pp-waves, our first  purpose in this work is to describe all non-flat left-invariant pp-wave metrics on four-dimensional Lie groups in Theorem~\ref{th:main}. With this description at hand, we finally determine all non-trivial four-dimensional Lie groups with harmonic curvature. 
A remarkable consequence of Theorem~\ref{th:main2}	is the existence of a single family of four-dimensional Lorentzian Lie groups with harmonic curvature whose Ricci tensor is not parallel and which  are neither locally conformally flat nor plane waves.

\subsection{pp-wave and plane wave Lorentzian manifolds}\label{pp and plane wave}
A  Lorentzian manifold $(M,g)$ admitting a parallel degenerate line field $\mathcal{U}$ is said to be a \emph{pp-wave} if the parallel line field is locally generated by a parallel null vector field $U$  and $(M,g)$ is transversally flat, i.e., $R(X,Y)=0$ for all $X,Y\in\mathcal{U}^\perp$. Any such Lorentzian four-manifold admits local coordinates $(x^+,x^-,x^1,x^2)$ so that 
 $$
 g=2 dx^+\circ dx^-+ H(x^+,x^1,x^2) dx^+\circ dx^+ +dx^1\circ dx^1+dx^2\circ dx^2\,,
 $$
where the degenerate parallel line field is generated by $U=\partial_{x^-}$.

The Ricci tensor of a pp-wave is of rank-one with $\rho({\partial x^+},{\partial x^+})=-\frac{1}{2}\Delta_xH$, where $\Delta_x=\partial^2_{x^1 x^1}+\partial^2_{x^2 x^2}$ is the space Laplacian. Hence the metric is Einstein (indeed Ricci-flat) if and only if the space Laplacian of $H$ vanishes.  As a general fact a four-dimensional pp-wave has harmonic curvature if and only if the defining function $H$ satisfies $\Delta_xH=\varphi(x^+)$ for some function $\varphi(x^+)$ (see, for example \cite{LN}). 

A pp-wave $(M,g,\mathcal{U})$  is a \emph{plane wave} if the covariant derivative of the curvature tensor satisfies $\nabla_XR~=~0$ for all $X\in\mathcal{U}^\perp$. In this case the local coordinates above can be specialized so that the function $H$ is quadratic on the space coordinates and, moreover, one may assume $H(x^+,x^1,x^2)=a_{ij}(x^+) x^ix^j$, from where it immediately follows that plane waves have harmonic curvature. It is clear that any locally symmetric pp-wave is a plane wave. Moreover pp-waves reduce to plane waves under some geometrical conditions as, for instance, in the locally conformally flat situation.

\subsection{Homogeneous setting}\label{homogeneous pp}

Homogeneous pp-waves have been broadly investigated (see, for example, \cite{GL} and references therein). They reduce to plane waves in many cases, although the existence of pp-waves is allowed for instance in the reducible case. This is partially due to the three-dimensional situation.

\subsubsection{Three-dimensional homogeneous pp-waves}\label{3D homogeneous pp}
Locally homogeneous pp-waves in dimension three have been described in \cite{GGN}. In local coordinates $(x^+,x^-,x)$ they are given by 
$g=2 dx^+\circ dx^-+ H(x^+,x) dx^+\circ dx^+ +dx\circ dx$, where the function $H(x^+,x)$ corresponds to one of the following:
\begin{itemize}
	\item[$(\mathcal{N}_b)$:] $H(x^+,x)=-2b^{-2} e^{b x}$, where $b\neq 0$, or
	\item[$(\mathcal{P}_c)$:] $H(x^+,x)=-\alpha(x^+)x^2$, where $\alpha'=c\alpha^{3/2}$, $c\neq 0$ and $\alpha>0$, or
	\item[$(\mathcal{CW})$:] $H(x^+,x)=-2 x^2$.
\end{itemize}
Metrics in classes $(\mathcal{P}_c)$ and $(\mathcal{CW})$ are plane waves, the latter corresponding to the three-dimensional Cahen-Wallach symmetric space. In contrast, metrics $(\mathcal{N}_b)$ are $pp$-waves but not plane waves.
All these metrics are realized as left-invariant metrics on three-dimensional Lie groups. Indeed, a straightforward calculation shows that a non-flat three-dimensional Lorentzian Lie group is a pp-wave if and only if it is isomorphically homothetic to the left-invariant metric determined by one of the following Lorentzian Lie algebras (where we follow the notation in \cite{ABDO}):
\begin{enumerate}
\item[(1)]The Lie algebra $\mathfrak{e}(1,1)$ with a Lorentzian metric described by
$[u_1,u_3]=\lambda u_1+ u_2$,\,\, $[u_2,u_3]=-\lambda u_2$, for any $\lambda\in\mathbb{R}$, $\lambda\neq 0$.	
\item[(2)] The Lie algebra $\mathfrak{e}(1,1)$ with a Lorentzian metric described by
$[u_1,u_2]=u_1$,\,\, $[u_2,u_3]=u_3$.
\item[(3)] The Lie algebra $\mathfrak{r}_{3,\lambda}$ and Lorentzian metric described by
$[u_1,u_2]= u_2$,\,\, $[u_1,u_3]=\gamma u_2+\lambda u_3$, for any $\gamma,\lambda\in\mathbb{R}$ with $\lambda(\lambda^2-1)\neq 0$.
\item[(4)] The Lie algebra $\mathfrak{r}_{3,0}=\mathfrak{aff}(\mathbb{R})\times\mathbb{R}$ with a Lorentzian metric described by
$[u_1,u_3]=\gamma u_2+ u_3$, for any $\gamma\in\mathbb{R}$,
\end{enumerate}
where $\{ u_i\}$ is a pseudo-orthonormal basis with  $\langle u_1,u_2\rangle=\langle u_3,u_3\rangle=1$.

Left-invariant metrics in class (1) are pp-waves corresponding to the local models $(\mathcal{N}_b)$, while those in (2) and (3) are non-symmetric plane waves locally modeled in $(\mathcal{P}_c)$. Finally left-invariant metrics in class (4) are Cahen-Wallach symmetric spaces corresponding to the model $(\mathcal{CW})$.
Since the connected and simply connected Lie group $E(1,1)$ admits cocompact subgroups (see \cite{Bock}), one has that the pp-wave structure (1) descends to the quotient solvmanifold, thus providing compact homogeneous pp-waves which are not plane waves.


\subsubsection{Four-dimensional homogeneous plane waves}\label{4D homogeneous plane wave}
A complete description of four-dimensional homogeneous pp-waves is not yet available, although the plane wave case is well understood.
Homogeneous plane waves in dimension four are described in terms of a $2\times 2$ skew-symmetric matrix $F$ and a $2\times 2$ symmetric matrix $A_0$ so that the defining function $H(x^+,x^1,x^2)$ takes the form
 $H=\vec{x}^T\, A(x^+) \,\vec{x}$, where $\vec{x}=(x^1,x^2)$ and the matrix $A(x^+)$ is given by (see \cite{BO03})
 $$ 
 (i)\,\,A(x^+)=e^{x^+\,F}A_0 e^{-x^+\, F},
 \quad\text{or}\quad
 (ii)\,\,A(x^+)=\frac{1}{(x^+)^2}e^{\log(x^+)F}A_0 e^{-\log(x^+) F}\,.
 $$
 Moreover, the plane wave metric is Ricci-flat if and only if $A_0$ is trace-free and locally conformally flat if $A_0$ is a multiple of the identity.
 
 Note that four-dimensional homogeneous plane waves of types (i) and (ii) are different representatives of the same conformal class, since for any homogeneous  plane wave $g=2 dx^+\circ dx^-+ H(x^+,\vec{x}) dx^+\circ dx^+ +\sum_i dx^i\circ dx^i$ of type (ii) with $H(x^+,\vec{x})=\vec{x}^T\, A(x^+) \,\vec{x}$ determined by a symmetric matrix $A_0$, the change of variables 
 $x^+\mapsto e^{x^+}$, \,
 $x^-\mapsto x^--\frac{1}{4}\sum (x^i)^2$,\,
 $x^i\mapsto e^{\frac{1}{2}x^+}x^i$,\,
 transforms the metric into a conformal one of the form
 $$
\widetilde{g}=e^{x^+}\left(  2 dx^+\circ dx^-+ \widetilde{H}(x^+,\vec{x}) dx^+\circ dx^+ +\sum_i dx^i\circ dx^i \right)
$$
 where $\widetilde{H}(x^+,\vec{x})=\vec{x}^T\, \widetilde{A}(x^+) \,\vec{x}$ corresponds to a type (i) plane wave determined by a symmetric matrix $\widetilde{A}_0=A_0-\frac{1}{4}\operatorname{id}$ (see \cite{HS}).

 A straightforward calculation shows that four-dimensional homogeneous plane waves in class (i) have parallel Ricci tensor, while metrics in class (ii) have parallel Ricci tensor if and only if they are Ricci-flat. 
 Furthermore, non-flat locally symmetric plane waves correspond to metrics in class (i) determined by a matrix $A(x^+)$ with constant coefficients. Note that  while non-symmetric three-dimensional homogeneous plane waves are of type (ii), it follows from the subsequent analysis that there are also plenty of non-symmetric left-invariant plane waves of type (i) in dimension four.

 \subsection{Harmonic Weyl tensor}
 The \emph{Schouten tensor} of an $n$-dimensional pseudo-Riemannian manifold $(M,g)$, given by $S=\rho-\frac{\tau}{2(n-1)}g$, determines the \emph{Cotton tensor} defined by $C(X,Y,Z)=(\nabla_XS)(Y,Z)-(\nabla_YS)(X,Z)$. The Cotton tensor is conformally invariant in dimension three, characterizing locally conformally flat three-dimensional manifolds as those with $C=0$. In higher dimensions it is not conformally invariant but $C=-\frac{ n-2}{n-3}\operatorname{div}W$, thus showing that the Cotton tensor vanishes if and only if the Weyl tensor is harmonic. Hence locally conformally flat manifolds, as well as those with parallel Ricci tensor have harmonic Weyl tensor. 
 
 The curvature tensor is harmonic if and only if the Ricci tensor is Codazzi, i.e., $(\nabla_X\rho)(Y,Z)-(\nabla_Y\rho)(X,Z)=0$, thus showing that harmonicity of the Weyl tensor and the curvature tensor are equivalent if the scalar curvature is constant (see \cite{Gr}).
 Lorentzian Lie groups with harmonic curvature have already been considered in \cite{AB}, where the Authors focus on some possibilities for the Jordan normal forms of the Ricci operator, without considering all the possible cases as shown in the following examples.

 \subsubsection{Cotton-flat three-dimensional homogeneous manifolds}\label{section-3D-Cotton flat}
 Locally conformally flat homogeneous three-dimensional manifolds with diagonalizable Ricci operator are locally symmetric.
  Three-dimensional locally conformally flat homogeneous manifolds which are neither locally symmetric nor plane waves are locally isomorphically homothetic to the Lorentzian Lie groups determined by (see \cite{HT}):
 \begin{enumerate}
 	\item[(1)]The Lie algebra $\mathfrak{sl}(2,\mathbb{R})$ with a Lorentzian metric described by
 	$$
 	[e_1,e_2]= e_3-\sqrt{3} e_2,\quad
 	[e_1,e_3]= e_2+\sqrt{3} e_3,\quad
 	[e_2,e_3]= 2e_1,
 	$$	
 	\item[(2)] the Lie algebra $\mathfrak{r}_{3,3}$ with a Lorentzian metric described by
 	$$
 	[e_1,e_2]=(\tfrac{1}{2}+\varepsilon\gamma)e_2+\gamma e_3,\quad
 	[e_1,e_3]=(\varepsilon-\gamma) e_2+(\tfrac{3}{2}-\varepsilon\gamma) e_3,\quad
 	\varepsilon^2=1,
 	$$
 \end{enumerate}
 	for any $\gamma\in\mathbb{R}$, $\gamma\neq\frac{\varepsilon}{2}$,  
 where $\{ e_i\}$ is an orthonormal basis with $e_3$ timelike.

 	The Ricci operator of any left-invariant metric (1) has real and complex eigenvalues $\{ -8, 4(1\pm\sqrt{3}\sqrt{-1})\}$, while the Ricci operator of metrics (2) has a single eigenvalue $-2$, which is a double root of the minimal polynomial. Left-invariant metrics determined by (1) have also been considered in \cite[Proposition 4.1]{AB}. However the case (2), having a single Ricci eigenvalue, does not correspond to any of the situations investigated in \cite{AB}.

Any three-dimensional locally conformally flat manifold $(N,h)$ trivally extends to a four-dimensional Cotton-flat product $(N\times\mathbb{R},h\oplus\pm dt^2)$ which is not necessarily locally conformally flat. 
Hence, we say that a four-dimensional Cotton-flat Lorentzian manifold is \emph{non-trivial} if and only if the Ricci tensor is not parallel and the manifold is not locally conformally flat, nor a product, nor a plane wave.

\subsection{Summary of results}
The purpose of this work is twofold. Firstly we determine all non-flat four-dimensional left-invariant metrics on Lie groups which are pp-waves.  Secondly we analyze the harmonicity of the Weyl tensor, obtaining a full description of four-dimensional Lorentzian Lie groups with harmonic curvature. 

Our first result shows that pp-wave Lorentzian Lie groups correspond to plane waves or to the Lie group $E(1,1)\times\mathbb{R}$ with the product metric of the three-dimensional homogeneous pp-wave discussed in Section~\ref{3D homogeneous pp} and the real line. This result is summarized in the following

\begin{theorem}\label{th:main}
A four-dimensional  non-flat left-invariant pp-wave Lie group
is isomorphically homothetic to one of the following:
\begin{itemize}
	\item[(1)]   
	$[u_1,u_4] = \gamma_1 u_1 +  u_2-\gamma_3 u_3$,
	$[u_2,u_4] = -\gamma_1 u_2$,
	$[u_3,u_4] = \gamma_3 u_2$, $\gamma_1\neq 0$.
		
	\smallskip
	
	\item[(2)]
	$[u_2,u_4] = -(\gamma_2-1)u_3$,
	$[u_3,u_4] = (\gamma_2+1)u_1$.

	\smallskip
	
	\item[(3)]
	$[v_1,v_4] = \gamma_1 v_1 -\gamma_2 v_2 + \gamma_3 v_3$,
	$[v_2,v_4] = \gamma_2 v_1 +\gamma_4 v_2 + \gamma_5 v_3$,
	$[v_3,v_4] = \gamma_8 v_3$.
	
	\smallskip
	
	\item[(4)]
	$[v_1,v_2]=v_3$, 
	$[v_1,v_4]= \gamma_1 v_1 -  \gamma_2 v_2 + \gamma_3 v_3$,
	$[v_2,v_4]=  \gamma_2 v_1 + \gamma_4 v_2 + \gamma_5 v_3$,
	$[v_3,v_4]=(\gamma_1+\gamma_4) v_3$.
	
	\smallskip
	
	\item[(5)]
	$[v_1,v_3]=v_2$,
	$[v_2,v_3]=-v_1$,
	$[v_1,v_4]=\gamma_1 v_1 + \gamma_2 v_2$,
	$[v_2,v_4]=-\gamma_2 v_1+\gamma_1 v_2$, $\gamma_1\neq 0$.
\end{itemize}

\smallskip

\noindent
Moreover, left-invariant metrics in cases~{\rm (2)}, {\rm (3), (4)} and {\rm (5)}  are plane waves.
Here $\{u_i\}$ and $\{v_i\}$  are   pseudo-orthonormal bases,   with $\langle u_1,u_2\rangle=\langle u_3,u_3\rangle=\langle u_4,u_4\rangle=1$
and 
$\langle v_1,v_1\rangle=\langle v_2,v_2\rangle=\langle v_3,v_4\rangle=1$, respectively.	
\end{theorem}

Left-invariant metrics in Theorem~\ref{th:main}--(1) are the only left-invariant pp-waves which are not plane waves. The underlying Lie algebra is the product $\mathfrak{e}(1,1)\times\mathbb{R}$ and the metric is isometric to the Riemannian direct extension of the three-dimensional homogeneous pp-wave in $E(1,1)$ discussed in Section~\ref{3D homogeneous pp} (see Remark~\ref{re:2.3}).

All Lie algebras in Theorem~\ref{th:main} are solvable. Moreover the Lie algebra in case (5) is $\mathfrak{aff}(\mathbb{C})$, while Lie algebras in case (2) are $\mathfrak{h}_3\times\mathbb{R}$ (resp., $\mathfrak{n}_4$) if $\gamma_2=\pm 1$ (resp., $\gamma_2\neq\pm 1$). Left-invariant metrics corresponding to cases (3) and (4) are realized on all semi-direct extensions of the Abelian and the Heisenberg algebras except the Lie algebra $\mathfrak{r}_{4,1,1}$ and $\mathbb{R}^4$, where any left-invariant metric is of constant sectional curvature (see \cite{Nomizu}).

Any left-invariant Lorentz metric on the four-dimensional Abelian Lie group $\mathbb{R}^4$ is flat, as well as any left-invariant Lorentz metric with nilpotent Ricci operator on a Lie group of  type~$\mathfrak{S}$ corresponding to the Lie algebra $\mathfrak{r}_{4,1,1}$ (see \cite{Nomizu}). Besides those cases, only the product affine group corresponding to the Lie algebra $\mathfrak{aff}(\mathbb{R})\times\mathfrak{aff}(\mathbb{R})$ does not admit any pp-wave left-invariant metric (see Section~\ref{se:EE}). Hence one has that
\begin{quote}\emph{
\hspace*{0.4cm}Any non-Abelian four-dimensional Lie group admits non-flat left-invariant Lorentzian metrics which are plane waves, with the exception of the Lie groups associated to the Lie algebras $\mathfrak{r}_{4,1,1}$ and $\mathfrak{aff}(\mathbb{R})\times\mathfrak{aff}(\mathbb{R})$.}
\end{quote}

\medskip 
	
In sharp contrast with plane waves, the existence of non-trivial Lorentzian Lie groups with harmonic curvature is a very restrictive condition. In addition to products with a three-dimensional locally conformally flat factor, plane waves,  locally conformally flat cases and metrics with parallel Ricci tensor (see \cite{AY, CZ2, CZ, CGGV}), four-dimensional Lorentzian Lie groups with harmonic curvature reduce to a family of examples with underlying Lie algebra $\mathfrak{r}_{4,\mu,-\mu}$ with $\mu=\frac{1}{4}(1\pm\sqrt{5})$ as follows. 

\begin{theorem}\label{th:main2}
A non-trivial four-dimensional  Lorentzian Lie group with harmonic curvature is isomorphically homothetic to 
$$
[u_1,u_4]=-\tfrac{1}{4}\eta_2(1+\varepsilon_1\sqrt{5})u_1+\varepsilon_2 u_2,\quad
[u_2,u_4]=\tfrac{1}{4}\eta_2(1+\varepsilon_1\sqrt{5})u_2,\quad
[u_3,u_4]=\eta_2 u_3,
$$
where $\eta_2\neq 0$, $\varepsilon_1^2=\varepsilon_2^2=1$,  and $\{ u_i\}$ is a pseudo-orthonormal basis with $\langle u_1,u_2\rangle=\langle u_3,u_3\rangle=\langle u_4,u_4\rangle=1$.
\end{theorem}

The Ricci operator of any metric in Theorem~\ref{th:main2} has eigenvalues $\{0,0,-\eta_2^2,-\eta_2^2\}$, the eigenvalue $0$ being a double root of the minimal polynomial. Due to the Jordan normal form of the Ricci operator, the metric above  corresponds to the situation in \cite[Theorem 5.1]{AB}.

As a consecuence of Theorem~\ref{th:main} and Theorem~\ref{th:main2}, one has that 
\begin{quote}
\emph{A  four-dimensional pp-wave Lie group is a plane wave if and only if the curvature is harmonic.}
\end{quote}

\subsection{Structure of the paper} 
Connected and simply connected non-solvable Lie groups in dimension four reduce to $SU(2)\times\mathbb{R}$ and $\widetilde{SL}(2,\mathbb{R})\times\mathbb{R}$. In the solvable case, they are isomorphic to semi-direct extensions of the three-dimensional unimodular Lie groups (see, for example,  \cite{ABDO}), thus corresponding to $\mathbb{R}^3\rtimes\mathbb{R}$, $\mathcal{H}^3\rtimes\mathbb{R}$, $E(1,1)\rtimes\mathbb{R}$ and $\widetilde{E}(2)\rtimes\mathbb{R}$, where $\mathbb{R}^3$ is the Abelian Lie group, $\mathcal{H}^3$ is the Heisenberg group, and $E(1,1)$ (resp., $\widetilde{E}(2)$) is the Poincaré (resp., Euclidean) Lie group.
In order to discuss left-invariant metrics, we work at the Lie algebra level  considering the different possibilities for Lorentzian metrics $\ip$ on Lie algebras $\mathfrak{g}=\mathfrak{k}\oplus\mathbb{R}$. 

Since any pp-wave has two-step nilpotent Ricci operator and has Petrov type N, we analyze these properties for all left-invariant Lorentz metrics on four-dimensional Lie groups. We then study the Cotton tensor towards a description of Lorentzian Lie groups with harmonic curvature. 
The analysis splits into three general situations depending on whether the restriction of $\ip$ to the unimodular ideal $\mathfrak{k}$ is Riemannian, Lorentzian or degenerate.

The discussion of Lorentzian metrics restricting to Riemannian metrics on $\mathfrak{k}$ follows from  the work of Milnor \cite{Milnor}, while the case when the restriction to $\mathfrak{k}$ is Lorentzian is based on the work of Rahmani \cite{Rahmani}. If $\ip \mid_{\mathfrak{k}\times\mathfrak{k}}$ is degenerate, then we consider the different possibilities on the dimension of the derived algebra $\mathfrak{k}'=[\mathfrak{k},\mathfrak{k}]$, inspired by the discussion in \cite{CC} (see also \cite{Rosalia}).
Left-invariant Lorentzian metrics on $\mathbb{R}^3\rtimes\mathbb{R}$ are considered in Section~\ref{se:R3}. In Section~\ref{se:H} we analyze left-invariant metrics on $\mathcal{H}^3\rtimes\mathbb{R}$. Semi-direct extensions of the Euclidean and Poincaré Lie groups are considered in Section~\ref{se:EE}, and the non-solvable case corresponding to $SU(2)\times\mathbb{R}$ and $\widetilde{SL}(2,\mathbb{R})\times\mathbb{R}$ is treated in Section~\ref{se:SLSU}.

Finally we emphasize that for any of the Lorentzian left-invariant metrics the properties $\operatorname{Ric}^2=0$ or $\operatorname{div}W=0$ amounts to solve some systems of polynomial equations in the structure constants. Since the solutions of a system of polynomial equations depend on the ideal generated by the polynomials rather than on the polynomials themselves, we make use of the theory of Gröbner bases (see, for example,  \cite{Cox}) to find simple polynomials in the corresponding ideals. As a technical fact, all the Gröbner bases calculations in this paper have been done with  {\sc Mathematica 12}, and computed using the lexicographical order with respect to the variables, which we specify in each case. 
We provide a complementary file with the calculations of the Gröbner bases, which is available at \texttt{https://doi.org/10.5281/zenodo.18548828}.
Observe that although $\operatorname{Ric}^2=0$ implies the vanishing of the scalar curvature $\tau$, the polynomial corresponding to  $\tau=0$ does not necessarily belongs to the ideal generated by the polynomials determined by the two-step nilpotency of the Ricci operator. Hence we explicitly use this polynomial determined by $\tau$ in the computation of the Gröbner bases when studying the property $\operatorname{Ric}^2=0$.

\section{Semi-direct extensions of the Abelian Lie group}\label{se:R3}

We consider separately the cases when the restriction of the metric to the three-dimensional Abelian ideal $\mathbb{R}^3$ is Riemannian (Section~\ref{R3-se: R3 Riemann}), Lorentzian (Section~\ref{R3-se:Lorentz}) or degenerate (Section~\ref{R3-se:degenerate}).  
The proof of the following theorem follows from the analysis below.
\begin{theorem}\label{R3-th: pp-waves}
	A four-dimensional  non-flat semi-direct extension of the Abelian Lie group   which is a pp-wave
	is isomorphically homothetic to one of the following:
	\begin{itemize}
		\item[(1)]   
		$[u_1,u_4] = \gamma_1 u_1 +  u_2-\gamma_3 u_3$,
		$[u_2,u_4] = -\gamma_1 u_2$,
		$[u_3,u_4] = \gamma_3 u_2$, $\gamma_1\neq 0$.

		\smallskip
		
		\item[(2)]
		$[u_2,u_4] = -(\gamma_2-1)u_3$,
		$[u_3,u_4] = (\gamma_2+1)u_1$.

		\smallskip
		
		\item[(3)]
		$[v_1,v_4] = \gamma_1 v_1 -\gamma_2 v_2 + \gamma_3 v_3$,
		$[v_2,v_4] = \gamma_2 v_1 +\gamma_4 v_2 + \gamma_5 v_3$,
		$[v_3,v_4] = \gamma_8 v_3$,
		
	\end{itemize}
	
	\smallskip
	
	\noindent
		where $\{u_i\}$ and $\{v_i\}$  are   pseudo-orthonormal bases,   with $\langle u_1,u_2\rangle=\langle u_3,u_3\rangle=\langle u_4,u_4\rangle=1$
	and 
	$\langle v_1,v_1\rangle=\langle v_2,v_2\rangle=\langle v_3,v_4\rangle=1$, respectively.
	
 Moreover, any non-trivial Cotton-flat semi-direct extension of the Abelian group is isomorphically homothetic to
	 	$$
	 	[u_1,u_4]=-\tfrac{1}{4}\eta_2(1+\varepsilon_1\sqrt{5})u_1+\varepsilon_2 u_2,\quad
	 	[u_2,u_4]=\tfrac{1}{4}\eta_2(1+\varepsilon_1\sqrt{5})u_2,\quad
	 	[u_3,u_4]=\eta_2 u_3,
	 	$$
	 	where $\eta_2\neq 0$, $\varepsilon_1^2=\varepsilon_2^2=1$,  and $\{ u_i\}$ is a pseudo-orthonormal basis as above.
\end{theorem}

\begin{remark}\rm
(i) Left-invariant metrics in case~(1) above, which are not plane waves, are never either Einstein nor locally conformally flat.

(ii) Left-invariant metrics in Theorem~\ref{R3-th: pp-waves}--(2) are Einstein if and only if $\gamma_2=0$, and they are locally conformally flat if and only if $\gamma_2=1$. Moreover, since the Ricci tensor is parallel, they are plane waves of type (i) whenever $\gamma_2\neq 0$.
	
(iii) For left-invariant metrics in case~(3) the Einstein condition is satisfied if and only if
$
\gamma_1^2  +  \gamma_4^2 - 
(\gamma_1 + \gamma_4) \gamma_8 = 0
$.  Moreover, these metrics are locally conformally flat if and only if either $\gamma_1=\gamma_4$ or $\gamma_1\neq \gamma_4$, $\gamma_2=0$ and $\gamma_8=\gamma_1+\gamma_4$. The Ricci tensor is parallel if and only if the metric is Einstein or $\gamma_8=0$. Hence, in the non-Einstein case, the plane waves are of type (i) if $\gamma_8=0$, and of type (ii) otherwise.
\end{remark}

\begin{remark}\rm \label{re:2.3}
A straightforward calculation shows that $\xi=\frac{\gamma_3}{\gamma_1} u_2+u_3$ determines a left-invariant unit spacelike parallel vector field on the Lie group corresponding to Theorem~\ref{R3-th: pp-waves}--(1). 
Moreover, considering the new basis
$\bar u_1= u_1-\frac{\gamma_3^2}{2\gamma_1^2}u_2-\frac{\gamma_3}{\gamma_1}u_3$,\,\,
$\bar u_2= u_2$,\,\,
$\bar u_3= u_4$,\,\,
$\bar u_4=\xi$,
the inner product is preserved and the Lie bracket transforms as
$[\bar u_1,\bar u_3]=\gamma_1 \bar u_1+\bar u_2$  and
$[\bar u_2,\bar u_3]=-\gamma_1\bar u_2$. 
This shows that the underlying Lie group is the product $E(1,1)\times\mathbb{R}$ and the Lorentzian metric is the product of the metric corresponding to (1) in Section~\ref{3D homogeneous pp} and the real line $\mathbb{R}\xi$.

For a given basis $\{ x_i\}$ of a Lie algebra $\mathfrak{g}$, let $\Lambda^2\mathfrak{g}$ be the space of bivectors with basis
$\left\{ x_1\wedge x_2,x_1\wedge x_3,x_1\wedge x_4,x_2\wedge x_3,x_2\wedge x_4,x_3\wedge x_4  \right\}$. We  consider the induced inner product given by 
$\langle\!\langle x_i\wedge x_j,x_k\wedge x_\ell\rangle\!\rangle=\langle x_i,x_k\rangle\langle x_j,x_\ell\rangle -\langle x_i,x_\ell\rangle\langle x_j,x_k\rangle$ so that the curvature tensor acts on the space of bivectors as
$\langle\!\langle \mathcal{R}(x_i\wedge x_j),x_k\wedge x_\ell\rangle\!\rangle=R(x_i,x_j,x_k,x_\ell)$.
Then it follows that the only non-zero component of $\mathcal{R}$ is given by
$\langle\!\langle\mathcal{R}(u_1\wedge u_4),u_1\wedge u_4\rangle\!\rangle=-2\gamma_1$, from where  $\mathcal{R}(u_1\wedge u_4)=-2\gamma_1u_2\wedge u_4$.	
\end{remark}

\subsection{Semi-direct extensions with Riemannian Lie group $\pmb{\mathbb{R}}^{\bm{3}}$} \label{R3-se: R3 Riemann}

Let $\mathfrak{g}=\algsdR$  be a semi-direct extension of the Abelian Lie algebra $\mathbb{R}^3$. If  $\ip$ is a Lorentzian inner product on $\mathfrak{g}$ whose restriction to $\mathbb{R}^3$ is of Riemannian signature then    there exists  an orthonormal basis $\{e_1,e_2,e_3,e_4\}$ of $\mathfrak{g}$, with $e_4$ timelike,   so that the structure of the metric Lie algebra is given by
\begin{equation}\label{R3-riem}
\begin{array}{l}
[e_1,e_4]=\eta_1 e_1-\gamma_1 e_2 - \gamma_2 e_3, \qquad
[e_2,e_4]=\gamma_1 e_1+\eta_2 e_2 - \gamma_3 e_3,
\\
\noalign{\medskip}
[e_3,e_4]=\gamma_2 e_1+\gamma_3 e_2 + \eta_3 e_3, 
\end{array} 
\end{equation}
for certain $\eta_i,\gamma_i\in\mathbb{R}$. 
In this case  a direct calculation shows that 
\[
\begin{array}{ll}
	\Riccicomp{1}{1} =  (\eta_1+\eta_2+\eta_3)\eta_1,
	&
	\Riccicomp{2}{1} = \Riccicomp{1}{2}=  \gamma_1(\eta_1-\eta_2),

	\\ \noalign{\medskip}
	
	\Riccicomp{2}{2} =  (\eta_1+\eta_2+\eta_3)\eta_2,
	&
	  \Riccicomp{3}{1} = \Riccicomp{1}{3} = \gamma_2(\eta_1-\eta_3),

	\\ \noalign{\medskip}
	
	\Riccicomp{3}{3} =   (\eta_1+\eta_2+\eta_3)\eta_3,

	&
	
	\Riccicomp{3}{2} = \Riccicomp{2}{3} = \gamma_3(\eta_2-\eta_3),
	
	\\ \noalign{\medskip}

	\Riccicomp{4}{4} =  \eta_1^2+\eta_2^2+\eta_3^2,
	&
	
\end{array}
\]
and therefore the scalar curvature is given by
$\tau= 2( \eta_1^2+\eta_2^2+\eta_3^2+\eta_1\eta_2+\eta_1\eta_3+\eta_2\eta_3)$. Note that $\tau$ vanishes if and only if $\eta_1=\eta_2=\eta_3=0$, in which case the space is flat.

\begin{remark}\label{re: Cotton Riemanniano}\rm
Since the scalar curvature of any left-invariant metric is constant, the Cotton tensor is equivalently given by $\Cotton(X,Y,Z)=(\nabla_X\rho)(Y,Z)-(\nabla_Y\rho)(X,Z)$. We denote by $\Cotton_{ijk}=\Cotton(e_i,e_j,e_k)$ the components of the Cotton tensor on the basis of left-invariant vector fields determined by $\{ e_i\}$.
Next we show that any left-invariant metric described by \eqref{R3-riem} is Cotton-flat if and only if it is Einstein or locally symmetric.			
			We start calculating the following three-components of the Cotton tensor:
			\[
			\Cotton_{124} =-\gamma_1(\eta_1-\eta_2)^2, \quad
			\Cotton_{134} =-\gamma_2(\eta_1-\eta_3)^2, \quad
			\Cotton_{234} =-\gamma_3(\eta_2-\eta_3)^2.
			\] 
			Note that  the isometry   
			$(e_1, e_2, e_3, e_4)\mapsto(e_1,e_3,e_2,e_4)$ gives the correspondence  $(\eta_1,\eta_2,\eta_3,\gamma_1,\gamma_2,\gamma_3)\sim (\eta_1,\eta_3,\eta_2,\gamma_2,\gamma_1,-\gamma_3)$, while
			$( e_1, e_2, e_3, e_4)\mapsto(e_3,e_2,e_1,e_4)$ shows that $(\eta_1,\eta_2,\eta_3,\gamma_1,\gamma_2,\gamma_3)\sim (\eta_3,\eta_2,\eta_1,-\gamma_3,-\gamma_2,-\gamma_1)$. 
			Using these isometries and excluding the Einstein case ($\eta_1=\eta_2=\eta_3$), the vanishing of the components above leads to analyze the case $\eta_1=\eta_2\neq\eta_3$, $\gamma_2=\gamma_3=0$, 
			and the case $\eta_1\neq \eta_2\neq\eta_3$, $\eta_1\neq \eta_3$, $\gamma_1=\gamma_2=\gamma_3=0$.
			
			If $\eta_1=\eta_2\neq \eta_3$ and $\gamma_2=\gamma_3=0$, we get $\Cotton_{141} = \eta_1\eta_3(\eta_3-\eta_1)$. Thus, either $\eta_1$ or $\eta_3$ vanishes and the space is locally symmetric. Now, if all the $\eta_i$'s are different and  $\gamma_1=\gamma_2=\gamma_3=0$, then the Cotton tensor is determined by
			\[
			\Cotton_{141} = \eta_1 \left( \eta_2^2+\eta_3^2-\eta_1(\eta_2+\eta_3) \right),\quad
			\Cotton_{242} = \eta_2 \left( \eta_1^2+\eta_3^2-\eta_2(\eta_1+\eta_3) \right),	 
			\]
			from where we get 
			$
			\eta_1 \Cotton_{242} - \eta_2 \Cotton_{141} = \eta_1\eta_2(\eta_1+\eta_2+\eta_3)(\eta_1-\eta_2)  
			$. Since the $\eta_i$'s are different one easily checks that $\Cotton_{141}=\Cotton_{242}=0$ does not hold if   any of the  first three factors $\eta_1$, $\eta_2$ or  $\eta_1+\eta_2+\eta_3$ vanishes. Hence there are no non-trivial Cotton-flat left-invariant metrics given by \eqref{R3-riem}.
		\end{remark}

\subsection{Semi-direct extensions with Lorentzian Lie group  $\pmb{\mathbb{R}}^{\bm{3}}$} \label{R3-se:Lorentz}
Let $\mathfrak{g}=\algsdR$ be a semi-direct extension of the Abelian Lie algebra $\mathbb{R}^3$. If the restriction to $\mathbb{R}^3$ of a Lorentzian inner product $\ip$ on  $\mathfrak{g}$ is of Lorentzian 
signature then we   must consider the possible Jordan normal forms of the   self-adjoint part of the derivation which determines the semi-direct extension. Next  we analyze them  separately.

\subsubsection{\bf The self-adjoint part of the derivation  is diagonalizable}
In this case there exists an orthonormal basis
$\{ e_1,e_2,e_3,e_4\}$,  with $e_3$ timelike, 
so that  left-invariant metrics are described by
\begin{equation}\label{R3-lorentz-Ia}
\begin{array}{l}
[e_1,e_4]=\eta_1 e_1-\gamma_1 e_2 + \gamma_2 e_3, \qquad
[e_2,e_4]=\gamma_1 e_1+\eta_2 e_2 + \gamma_3 e_3,
\\
\noalign{\medskip}
[e_3,e_4]=\gamma_2 e_1+\gamma_3 e_2 + \eta_3 e_3, 
\end{array} 
\end{equation}
for certain $\eta_i,\gamma_i\in\mathbb{R}$.

\medskip

We proceed as in the previous case.
The Ricci operator is determined by 
\[
\begin{array}{ll}
	\Riccicomp{1}{1} = - (\eta_1+\eta_2+\eta_3)\eta_1,
	&
	 \Riccicomp{2}{1} =  \phantom{-}\Riccicomp{1}{2}= -\gamma_1(\eta_1-\eta_2),

	\\ \noalign{\medskip}
	
	\Riccicomp{2}{2} = - (\eta_1+\eta_2+\eta_3)\eta_2,
	&
	  \Riccicomp{3}{1} = - \Riccicomp{1}{3} = \phantom{-} \gamma_2(\eta_1-\eta_3),

	\\ \noalign{\medskip}
	
	\Riccicomp{3}{3} =  - (\eta_1+\eta_2+\eta_3)\eta_3,
	&
	\Riccicomp{2}{3} = - \Riccicomp{3}{2} = -\gamma_3(\eta_2-\eta_3),
	
	\\ \noalign{\medskip}
	
	\Riccicomp{4}{4} =  -(\eta_1^2+\eta_2^2+\eta_3^2),
	&
	
\end{array}
\]
so that the scalar curvature  
$\tau= -2( \eta_1^2+\eta_2^2+\eta_3^2+\eta_1\eta_2+\eta_1\eta_3+\eta_2\eta_3)$
vanishes if and only if $\eta_1=\eta_2=\eta_3=0$, in which case the space is flat.

	\begin{remark}\rm
		Although in this case we only have a basic isometry,
		$( e_1, e_2, e_3, e_4)\mapsto(e_2,e_1,e_3,e_4)$, which 
		shows that $(\eta_1,\eta_2,\eta_3,\gamma_1,\gamma_2,\gamma_3)\sim (\eta_2,\eta_1,\eta_3,-\gamma_1,\gamma_3,\gamma_2)$, 
		the Cotton-flat condition can be handled exactly as in Remark~\ref{re: Cotton Riemanniano}. Indeed, considering that isometry and the vanishing of the components
		\[
		\Cotton_{124} =\gamma_1(\eta_1-\eta_2)^2, \quad
		\Cotton_{134} =\gamma_2(\eta_1-\eta_3)^2, \quad
		\Cotton_{234} =\gamma_3(\eta_2-\eta_3)^2,
		\]
		we are led to three possibilities: the case $\eta_1=\eta_2\neq\eta_3$, $\gamma_2=\gamma_3=0$, 
		the case $\eta_1=\eta_3\neq\eta_2$, $\gamma_1=\gamma_3=0$, 
		and the case $\eta_1\neq \eta_2\neq\eta_3$, $\eta_1\neq \eta_3$, $\gamma_1=\gamma_2=\gamma_3=0$ (excluding the Einstein situation $\eta_1=\eta_2=\eta_3$).
		We proceed exactly as in Remark~\ref{re: Cotton Riemanniano} to show that any Cotton-flat metric determined by \eqref{R3-lorentz-Ia} is trivial. 
	\end{remark}

\subsubsection{\bf The self-adjoint part of the derivation  has complex eigenvalues}
There exists an orthonormal basis  $\{ e_1,e_2,e_3,e_4\}$  of $\mathfrak{g}=\algsdR$, with $e_3$ timelike, so that the corresponding left-invariant metrics are described by
\begin{equation}\label{R3-lorentz-Ib}
\begin{array}{l}
[e_1,e_4]=\eta e_1-\gamma_1 e_2 + \gamma_2 e_3,  \qquad\quad
[e_2,e_4]=\gamma_1 e_1+\delta e_2 + (\gamma_3-\nu) e_3,
\\
\noalign{\medskip}
[e_3,e_4]=\gamma_2 e_1+(\gamma_3+\nu) e_2 + \delta e_3, 
\end{array} 
\end{equation}
for certain $\eta,\delta,\gamma_i\in\mathbb{R}$ and $\nu\neq 0$.
A straightforward calculation shows that the Ricci operator  is determined by
\[
\begin{array}{ll}
	\Riccicomp{1}{1} =  -(\eta+2\delta)\eta,
	&
	\Riccicomp{2}{1} =\phantom{-} \Riccicomp{1}{2}=  \gamma_1(\delta-\eta)-\gamma_2\nu,
	
	\\ \noalign{\medskip}
	
	\Riccicomp{2}{2} = -\delta\eta -2\delta^2 -2\gamma_3\nu,
	&
	 \Riccicomp{3}{1}= - \Riccicomp{1}{3} =  \gamma_2(\eta-\delta)-\gamma_1\nu,
	
	\\ \noalign{\medskip}
	
	\Riccicomp{3}{3} = -\delta\eta -2\delta^2+2\gamma_3\nu,
	&
	\Riccicomp{3}{2}= - \Riccicomp{2}{3} = (\eta+2\delta)\nu,
	
	\\ \noalign{\medskip}
	
	\Riccicomp{4}{4} = 2\nu^2-\eta^2-2\delta^2,
	&
	
\end{array}
\]
so that   $\Riccicomp{4}{4}$ is an eigenvalue. Hence $\Riccicomp{4}{4}$ must vanish, so   
$\nu^2 = \frac{1}{2} (\eta^2+2\delta^2)$. 
In this setting the scalar curvature 
$\tau = -(\eta+2\delta)^2$ and therefore $\eta=-2\delta$.

At this point we consider the  $\Ricci^2$ operator, which must  vanish.
Taking into account that $\eta=-2\delta$,   we compute
$\RicciRiccicomp{2}{2} - \RicciRiccicomp{3}{3} = (\gamma_1^2+\gamma_2^2)(9\delta^2+\nu^2)$.
Since $\nu\neq 0$, it follows that   $\gamma_1=\gamma_2=0$ and a long but straightforward calculation shows that the Weyl tensor acting on the space of bivectors, $\mathcal{W}:\Lambda^2\rightarrow\Lambda^2$, satisfies
\[
\langle\!\langle	\mathcal{W}^2(e_1\wedge e_4),e_1\wedge e_4\rangle\!\rangle=\tfrac{4}{9}(3\delta^2+\nu^2)^2 \neq 0.
\]
Thus, we conclude that there is no left-invariant pp-wave metric in this case.
	
	\begin{remark}\rm
Any Cotton-flat left-invariant metric \eqref{R3-lorentz-Ib} is necessarily Einstein. To show it, we start calculating $\Cotton_{234}=-4\gamma_3\nu^2$, from where $\gamma_3=0$. Now, 
		\[
		\Cotton_{124} = \gamma_1 \left( (\delta-\eta)^2-\nu^2 \right)
		-2\gamma_2(\delta-\eta)\nu,
		\,\,\,
		\Cotton_{134} = \gamma_2 \left( (\delta-\eta)^2-\nu^2 \right)
		+2\gamma_1(\delta-\eta)\nu,
		\]
		imply
		$\gamma_1 \Cotton_{134} -\gamma_2 \Cotton_{124} = 
		2(\gamma_1^2+\gamma_2^2)(\delta-\eta)\nu$, 
		so that either $\delta=\eta$ or $\gamma_1=\gamma_2=0$. Moreover, if $\delta=\eta$ then $\Cotton_{124}=-\gamma_1\nu^2 $ and 
		$\Cotton_{134}= -\gamma_2\nu^2$. Therefore, in any case,  $\gamma_1=\gamma_2=0$ and under this assumption a straightforward calculation shows that the Cotton tensor is determined by
		\[
		\Cotton_{243}=(2\delta^2-\eta^2+2\delta\eta+2\nu^2)\nu,\quad
		\Cotton_{343}=-(\delta^2-\delta\eta-\nu^2)\eta.
		\]
		Finally, one easily checks that the vanishing of these components leads to $\eta=-2\delta=-\frac{2\nu}{\sqrt{3}}\ve$, with $\ve^2=1$, and the corresponding left-invariant metric is Einstein.
		
	\end{remark}

\subsubsection{\bf The self-adjoint part of the derivation   has a double root of the minimal polynomial}
There exists a pseudo-orthonormal basis 
$\{u_1,u_2,u_3,u_4\}$, with $\langle u_1,u_2\rangle=\langle u_3,u_3\rangle=\langle u_4,u_4\rangle=1$, so that
the corresponding left-invariant metrics are described by 
\begin{equation}\label{R3-lorentz-II}
\begin{array}{l}
[u_1,u_4]=(\eta_1+\gamma_1) u_1 + \ve  u_2 - \gamma_3 u_3, \,\,\,\,\,
[u_2,u_4]= (\eta_1-\gamma_1) u_2 -  \gamma_2 u_3,
\\
\noalign{\medskip}
[u_3,u_4]=\gamma_2 u_1+\gamma_3 u_2 + \eta_2 u_3 , 
\end{array} 
\end{equation}
where $\eta_i,\gamma_i\in\mathbb{R}$ and  $\varepsilon^2=1$.

\medskip

A straightforward calculation shows that the Ricci operator  is determined by
\[
\begin{array}{ll}
	\Riccicomp{1}{1} = \Riccicomp{2}{2} = -(2\eta_1+\eta_2)\eta_1,
	&
	\Riccicomp{2}{1} = -\ve (2\eta_1+\eta_2+2\gamma_1),

	\\ \noalign{\medskip}
	
	\Riccicomp{3}{3} = -(2\eta_1+\eta_2)\eta_2, 
	&
	\Riccicomp{3}{1}  = \Riccicomp{2}{3} = -\gamma_3(\eta_1-\eta_2) -\ve \gamma_2 ,

	\\ \noalign{\medskip}
	
	\Riccicomp{4}{4} =  -2\eta_1^2-\eta_2^2, 
	&
	\Riccicomp{3}{2}=  \Riccicomp{1}{3}=  -\gamma_2(\eta_1-\eta_2) .
			
\end{array}
\]
Therefore, the scalar curvature $\tau=-2(3\eta_1^2+\eta_2^2+2\eta_1\eta_2)$, so that necessarily $\eta_1=\eta_2=0$. 
Now  the only non-zero component of $\Ricci^2$ is $\RicciRiccicomp{2}{1}=\gamma_2^2$.
Hence, $\gamma_2=0$ and the left-invariant metric is given by
\begin{equation}\label{eq:R3 Ib pp-wave}
[u_1,u_4] = \gamma_1 u_1 + \ve u_2-\gamma_3 u_3,
\quad
[u_2,u_4] = -\gamma_1 u_2,
\quad
[u_3,u_4] = \gamma_3 u_2.
\end{equation}
Note that the isomorphic isometry $u_4\mapsto -u_4$ interchanges the sign of  $\ve$, $\gamma_1$ and $\gamma_3$. Thus, without loss of generality, we may assume $\ve=1$. 
Besides, a straightforward  calculation shows that the metric is flat if and only if $\gamma_1=0$. 

In the non-flat case, since $\nabla_{u_i}u_2=0$ for $i=1,2,3$ and $\nabla_{u_4}u_2= \gamma_1 u_2$, it follows that $u_2$ determines a left-invariant recurrent null vector field.
Moreover,   a direct calculation shows that the curvature tensor satisfies $R(x,y)=0$ for all $x,y\in u_2^\perp = \span\{u_2,u_3,u_4\}$. 
However, the covariant derivative
$ (\nabla_{u_4}R)(u_1,u_4,u_1)= -4\gamma_1^2u_4\neq 0$.
Hence, we conclude that  the underlying structure for left-invariant metrics determined by Equation~\eqref{eq:R3 Ib pp-wave} is a  pp-wave but not a plane wave (cf.~\cite{Leistner}),   corresponding to case~(1) in Theorem~\ref{R3-th: pp-waves}.

\medskip

%
%

 	\begin{remark}\rm
 		Considering the components of the Cotton tensor of any metric \eqref{R3-lorentz-II},
 		\[
 		\Cotton_{234}=\gamma_2 (\eta_1-\eta_2)^2,\quad
 		\Cotton_{142}= -\ve\gamma_2^2-(2\gamma_2\gamma_3-\eta_1\eta_2)(\eta_1-\eta_2),
 		\]
 		if the Cotton tensor vanishes, then
 		it follows that  $\gamma_2=0$ and $\eta_1\eta_2(\eta_1-\eta_2)=0$. Then one has  
 		$\Cotton_{134}= \gamma_3(\eta_1-\eta_2)^2$, so that we are led to three possibilities: the case $\eta_1=\eta_2$, the case $\gamma_3=\eta_2=0\neq\eta_1$, and the case $\gamma_3=\eta_1=0\neq\eta_2$.
 		
 		If $\eta_1=\eta_2$, the Cotton tensor is determined by $\Cotton_{141}=\ve(2\gamma_1+3\eta_1)(2\gamma_1+\eta_1)$ and the space is Einstein for $2\gamma_1+3\eta_1=0$ and locally conformally flat if $2\gamma_1+\eta_1=0$. 
 		
 		If $\gamma_3=\eta_2=0\neq\eta_1$, the Cotton tensor is given by $\Cotton_{141} = 2\ve(\gamma_1+\eta_1)(2\gamma_1+\eta_1)$. For $\gamma_1+\eta_1=0$   the metric is locally symmetric. Otherwise, if $2\gamma_1+\eta_1=0$, the left-invariant metric corresponds to a   product Lie algebra $\alg\times\mathbb{R}$, with $\alg=\operatorname{span}\{ u_1,u_2,u_4\}$,
 		and a direct calculation shows that   $\nabla_{u_i} u_3 = 0$, $i=1,\dots,4$. Hence, $u_3$ generates a left-invariant spacelike parallel vector field and  the Lorentzian Lie group splits as a Lorentzian product Lie group which corresponds to the Lie algebra $\mathfrak{r}_{3,3}\times\mathbb{R}$ with the metric on the three-dimensional factor given by (2) in Section~\ref{section-3D-Cotton flat}. 
 		
 		Finally, we analyze the case $\gamma_3=\eta_1=0\neq\eta_2$. In this last case, the Cotton tensor is determined by $\Cotton_{141}=\ve(4\gamma_1^2-\eta_2^2+2\gamma_1\eta_2)$. Taking 
 		$\gamma_1=-\frac{1}{4}(1\pm\sqrt{5})\eta_2\neq 0$, a straightforward calculation shows that the   left-invariant metric, given by
 		\[
 		[u_1,u_4]=\gamma_1 u_1 + \ve  u_2 , \quad
 		[u_2,u_4]= -\gamma_1 u_2 , \quad
 		[u_3,u_4]=  \eta_2 u_3 , 
 		\]
 		is non-trivial Cotton-flat, corresponding to the Lorentzian Lie group in Theorem~\ref{th:main2}.
 		
 	A direct calculation shows that the left-invariant vector field $U_2$ determined by $u_2$ is null and recurrent, which shows that the underlying structure is a Brinkmann wave which is not transversally flat on the non-unimodular Lie algebra $\mathfrak{r}_{4,\mu,-\mu}$ with $\mu=\frac{1}{4}(1\pm\sqrt{5})$.
 			Moreover, the Ricci operator has eigenvalues $0$ and $-\eta_2^2$, the former being a double root of the minimal polynomial. The Weyl curvature operator $\mathcal{W}:\Lambda^2\rightarrow\Lambda^2$ has eigenvalues $-\frac{1}{3}\eta_2^2$ with multiplicity two and $\frac{1}{6}\eta_2^2$ with multiplicity four, the latter one being a double root of the minimal polynomial. Hence of Petrov Type II.
 	\end{remark}

\subsubsection{\bf The self-adjoint part of the derivation   has a triple root}
There exists a  pseudo-orthonormal basis 
$\{u_1,u_2,u_3, u_4\}$  of    $\algsdR$, with $\langle u_1,u_2\rangle=\langle u_3,u_3\rangle=\langle u_4,u_4\rangle=1$, 
so that the left-invariant metrics are given by
\begin{equation}\label{R3-lie-lorentzIII} 
\begin{array}{l}
[u_1,u_4]=(\eta+\gamma_1) u_1  - \gamma_3 u_3,  \quad
[u_2,u_4]= (\eta-\gamma_1) u_2 -  (\gamma_2-1) u_3,
\\
\noalign{\medskip}
[u_3,u_4]=(\gamma_2+1) u_1+\gamma_3 u_2 + \eta u_3, 
\end{array} 
\end{equation}
where $\eta, \gamma_i\in\mathbb{R}$.

\medskip

In this case the Ricci operator  is determined by
\[
\begin{array}{ll}
	\Riccicomp{1}{1} = \Riccicomp{2}{2} = -3\eta^2+\gamma_3,  
	&
	\Riccicomp{1}{2}= 2\gamma_2,
	
	\\ \noalign{\medskip}
	
	\Riccicomp{3}{3} = -3\eta^2-2\gamma_3, 
	&
	\Riccicomp{1}{3} = \Riccicomp{3}{2} =  -3\eta+\gamma_1,	
	
	\\ \noalign{\medskip}
	
	\Riccicomp{4}{4} = -3\eta^2,

\end{array}
\]
so that $\tau=-12\eta^2$ which leads to $\eta=0$. Now a direct calculation shows that $\Ricci^2=0$ if and only if $\gamma_1=\gamma_3=0$, and the corresponding left-invariant metric is given by
\begin{equation}\label{eq: R3 III plane wave}
[u_2,u_4] = -(\gamma_2-1)u_3,
\quad
[u_3,u_4] = (\gamma_2+1)u_1. 
\end{equation}
A straightforward calculation shows that $u_1$ determines a left-invariant parallel null vector field and, moreover, the curvature tensor satisfies $R(x,y)=0$ for all $x,y\in u_1^\perp=\span\{u_1,u_3,u_4\}$, while the covariant derivatives $\nabla_x R=0$ for all $x\in u_1^\perp$. Hence it follows from~\cite{Leistner} that the underlying structure of left-invariant metrics \eqref{eq: R3 III plane wave} is a plane wave corresponding to  Theorem~\ref{R3-th: pp-waves}--(2).

	\begin{remark}\rm
The vanishing of the Cotton tensor of any left-invariant metric given by \eqref{R3-lie-lorentzIII} is determined by the vanishing of
		\[
		\Cotton_{143}=-3\gamma_3^2,\,
		\Cotton_{242}=\gamma_1(6\gamma_2-1)-(8\gamma_2-3)\eta,\,
		\Cotton_{243}=(\gamma_1-3\eta)(\gamma_1-\eta)-(5\gamma_2-2)\gamma_3.
		\] 
Hence one has three possibilities: 
		the case $\gamma_1=\gamma_3=\eta=0$, the case $\gamma_1=3\eta\neq 0$, $\gamma_2=\gamma_3=0$, and the case 
		$\gamma_1=\eta\neq 0$, $\gamma_2=1$, $\gamma_3=0$.
		Note that the first case corresponds to the plane waves given by~\eqref{eq: R3 III plane wave}. Moreover, the metric is Einstein in the second case and locally conformally flat in the third one.
	\end{remark}

\subsection{Semi-direct extensions with degenerate Lie group $\pmb{\mathbb{R}}^{\bm{3}}$} \label{R3-se:degenerate}

In this case we consider  a four-dimensional Lie algebra    $\mathfrak{g}=\algsdR$   with a Lorentzian inner product $\ip$ which restricts to a degenerate inner product on the subalgebra $\mathbb{R}^3$.    
There exists a pseudo-orthonormal basis $\{u_1,u_2,u_3,u_4\}$  of $\mathfrak{g}$, with $\langle u_1,u_1\rangle=\langle u_2,u_2\rangle=\langle u_3,u_4\rangle=1$, possibly after rotating the vectors $u_1$ and $u_2$, so that 
\begin{equation}\label{R3-lie-degenerate}
\begin{array}{l}
[u_1,u_4]=\gamma_1 u_1-\gamma_2 u_2 + \gamma_3 u_3, \qquad 
[u_2,u_4]=\gamma_2 u_1+\gamma_4 u_2 + \gamma_5 u_3,
\\
\noalign{\medskip}
[u_3,u_4]=\gamma_6 u_1+\gamma_7 u_2 + \gamma_8 u_3, 
\end{array} 
\end{equation}
for certain  $\gamma_i\in\mathbb{R}$.   
 
 \medskip

We start computing the Ricci operator, determined by
\[
\begin{array}{ll}
	\Riccicomp{1}{1} =  -\frac{\gamma_6^2}{2},\quad 	\Riccicomp{2}{2} = -\frac{\gamma_7^2}{2},
	&

	\Riccicomp{3}{1} = \Riccicomp{1}{4} = \frac{1}{2}(2\gamma_1\gamma_6-\gamma_2\gamma_7+\gamma_4\gamma_6),
	
	\\ \noalign{\medskip}
	
	\Riccicomp{3}{3} = \Riccicomp{4}{4} =  \frac{1}{2}(\gamma_6^2+\gamma_7^2),

	&
	
	\Riccicomp{3}{2} = \Riccicomp{2}{4} =
			\frac{1}{2}(\gamma_1\gamma_7+\gamma_2\gamma_6+2\gamma_4\gamma_7),		
	
	\\ \noalign{\medskip}
	
	\Riccicomp{2}{1} = \Riccicomp{1}{2}= -\frac{\gamma_6\gamma_7}{2},
	&
	\Riccicomp{3}{4} = -\gamma_1^2-\gamma_4^2+\gamma_1\gamma_8 
		-\gamma_3\gamma_6+\gamma_4\gamma_8-\gamma_5\gamma_7,

%
%
\end{array}
\]
from where we get $\tau= \frac{1}{2}(\gamma_6^2+\gamma_7^2)$. Hence, necessarily $\gamma_6=\gamma_7=0$, which ensures that $\Ricci^2=0$, and the associated left-invariant metric corresponds to
\begin{equation}\label{eq:R3 deg plane wave}
[u_1,u_4] = \gamma_1 u_1 -\gamma_2 u_2 + \gamma_3 u_3,
\quad
[u_2,u_4] = \gamma_2 u_1 +\gamma_4 u_2 + \gamma_5 u_3,
\quad
[u_3,u_4] = \gamma_8 u_3 .
\end{equation}
Now a straightforward calculation shows that  $u_3$ determines a left-invariant recurrent null vector field since 
$\nabla_{u_i}u_3=0$ for $i=1,2,3$ and $\nabla_{u_4}u_3= -\gamma_8 u_3$.  
Moreover, the curvature tensor satisfies $R(x,y)=0$ for all $x,y\in u_3^\perp = \span\{u_3,u_1,u_2\}$, while the covariant derivative satisfies 
$ \nabla_x R=0$ for all $x\in u_3^\perp$.
Hence, we conclude that the underlying structure of  left-invariant metrics determined by Equation~\eqref{eq:R3 deg plane wave} is a  plane wave (cf.~\cite{Leistner}),  corresponding to case~(3) in Theorem~\ref{R3-th: pp-waves}.

 	\begin{remark}\rm
 		The components $
 		\Cotton_{134}=-\tfrac{1}{2}(\gamma_6^2+\gamma_7^2)\gamma_6$ and
 		$\Cotton_{234}=-\tfrac{1}{2}(\gamma_6^2+\gamma_7^2)\gamma_7$ of the Cotton tensor of any metric \eqref{R3-lie-degenerate} imply $\gamma_6=\gamma_7=0$. Hence the Weyl tensor is harmonic if and only if the underlying structure is a plane wave  given by~\eqref{eq:R3 deg plane wave}.
 	\end{remark}

\section{Semi-direct extensions of the Heisenberg Lie group $\mathcal{H}^3$}\label{se:H}

Semi-direct extensions of the Heisenberg algebra include the nilpotent Lie algebras $\mathfrak{h}_3\times\mathbb{R}$ and $\mathfrak{n}_4$, which are also semi-direct extensions of the Abelian Lie algebra. Since these two cases have already been considered in the discussion of the previous section, we focus on semi-direct extensions of the Heisenberg group which are not almost Abelian (see \cite{KT} for a description of left-invariant metrics on the product Lie group $\mathcal{H}^3\times\mathbb{R}$).

We proceed as in the case of the Abelian Lie group, considering separately the cases when the induced metric on $\mathcal{H}^3$ is positive definite, Lorentzian or degenerate.  
In this case, the use of Gröbner bases facilitates obtaining the following result, which is derived from the analysis we carry out below.

\begin{theorem}\label{H3-th: pp-waves}
	A four-dimensional    non-flat semi-direct extension of the Heisenberg Lie group  which is a pp-wave not realized on $\algsdR$    is   isomorphically homothetic to  a plane wave
	\[
	\begin{array}{ll}
		[v_1,v_2]=    v_3 , &
		[v_1,v_4]= \gamma_1 v_1 -  \gamma_2 v_2 + \gamma_3 v_3 ,
		\\
		\noalign{\medskip}
		[v_2,v_4]=  \gamma_2 v_1 + \gamma_4 v_2 + \gamma_5 v_3 , &
		[v_3,v_4]=(\gamma_1+\gamma_4) v_3 ,
	\end{array}
	\]
	
	\smallskip
	
	\noindent
where $\{v_i\}$ is a pseudo-orthonormal basis   with $\langle v_1,v_1\rangle=\langle v_2,v_2\rangle=\langle v_3,v_4\rangle=1$.

Moreover, any Cotton-flat semi-direct extension of the Heisenberg group which is not realized on $\algsdR$ is locally conformally flat or a plane wave as above.
\end{theorem}

\begin{remark}\rm
	Left-invariant metrics given in the above theorem are Einstein if and only if $4\gamma_1\gamma_4+1=0$, and they are locally conformally flat if and only if 
	$(\gamma_1-\gamma_4)(2\gamma_2-1)=0$. Moreover, the Ricci tensor is parallel if and only if the metric is Einstein or $\gamma_1+\gamma_4=0$. In the non-Einstein case one has therefore that plane waves in Theorem~\ref{H3-th: pp-waves} are of type (i) if $\gamma_1+\gamma_4=0$, and of type (ii) otherwise.
	
	Considering the new basis
	$\bar v_1 = v_1$, $\bar v_2 = v_2$, $\bar v_3 = v_3$,
	$\bar v_4 =     \gamma_5 v_1 -   \gamma_3  v_2  +v_4$,
	the Lie bracket transforms into
	\[
	[\bar v_1,\bar v_2] =       \bar   v_3,
	\,\,
	[\bar v_1,\bar v_4] = \gamma_1 \bar v_1 -\gamma_2 \bar v_2,
	\,\,
	[\bar v_2, \bar v_4] =   \gamma_2 \bar v_1 + \gamma_4 \bar v_2 ,
	\,\,
	[\bar v_3, \bar v_4] = (\gamma_1+\gamma_4) \bar v_3,
	\]
	and a   direct  calculation shows that, when evaluating on the basis $\{\bar v_i\}$,
	\[ 
	\ad_{\bar v_4}=\left(
	\begin{array}{cc}
		A& 0 
		\\
		0 & \operatorname{tr}A
	\end{array}
	\right), 
	\quad \text{where}\,\, 
	A=\left(
	\begin{array}{cc}
		-\gamma_1 & -\gamma_2  
		\\
		\gamma_2 & -\gamma_4  
	\end{array}
	\right).
	\qquad\qquad
	\]
	Hence,  the Lie algebra corresponds to $ \mathfrak{h}_3\times\mathbb{R} $ or $\mathfrak{n}_4$ if and only if $\gamma_4=-\gamma_1$ and 
	$\gamma_1^2=\gamma_2^2$ (cf. \cite{ABDO}). 
	Indeed, setting $\gamma_2=\ve \gamma_1$, with $\ve^2=1$,   and considering the basis 
	\[
	\widetilde v_1 = -  \gamma_5     v_1 
	+  \gamma_3    v_2  -  v_4,\quad
	\widetilde v_2 =       v_3,\quad
	\widetilde v_3 = \ve v_1 - v_2,\quad
	\widetilde v_4 =  v_1,
	\]
	the Lie brackets transform into
	$[\widetilde v_1,\widetilde v_4] =   \ve \gamma_1 \widetilde v_3$ and 
	$[\widetilde v_3, \widetilde v_4] =  \widetilde v_2$.
	Therefore,   we conclude that only if  $\gamma_4=-\gamma_1$ and $\gamma_1^2 = \gamma_2^2$   the semi-direct extensions in Theorem~\ref{H3-th: pp-waves} are almost Abelian and therefore covered by Theorem~\ref{R3-th: pp-waves}.

\end{remark}

\subsection{Semi-direct extensions with Riemannian normal subgroup $\bm{\mathcal{H}^3$}} \label{H3-se: H3 Riemann}

In this section we consider   left-invariant  Lorentzian metrics which are obtained as extensions 
of the three-dimensional Riemannian Heisenberg Lie group $\mathcal{H}^3$. 
One may describe all such Lorentz extensions considering an orthonormal basis $\{ e_1,e_2,e_3,e_4\}$, with $e_4$ timelike, where the Lie brackets become
\begin{equation}\label{H3-riemannian}
\begin{array}{ll}
[e_1,e_2]=\lambda_3 e_3,  &
[e_1,e_4]=\gamma_1 e_1-\gamma_2 e_2 + \gamma_3 e_3,  
\\
\noalign{\medskip}
[e_2,e_4]=\gamma_2 e_1 + \gamma_4  e_2 + \gamma_5 e_3,  &
[e_3,e_4]=(\gamma_1 +\gamma_4)e_3,
\end{array}  
\end{equation}
where $\lambda_3\neq 0$ and $\gamma_i\in\mathbb{R}$.

 \medskip

We start considering an   auxiliary variable $\lambda'_3$ to indicate that $\lambda_3\neq 0$	by means of  $\lambda_3 \lambda'_3-1$.
In the polynomial ring $\mathbb{R}[\lambda'_3,\gamma_3,\gamma_5,\gamma_2,\gamma_1,\gamma_4,\lambda_3]$ we compute a Gröbner basis for the ideal $\mathcal{I}
= \langle \RicciRiccicomp{i}{j} \cup \{\tau, \lambda_3\lambda'_3-1 \} \rangle$  and   get $19$ polynomials, among which we find
$
\mathbf{g}_1 = \gamma_2^2$,
$\mathbf{g}_2 =(4\gamma_4^2+\lambda_3^2)(16\gamma_4^2+\lambda_3^2)\gamma_4^2$ and
$\mathbf{g}_3 = (\gamma_1+2\gamma_4)(2\gamma_1+\gamma_4)$.
%
%
%
%
%
%
Hence, since $\lambda_3\neq 0$,  necessarily $\gamma_1=\gamma_2=\gamma_4=0$ and the possible left-invariant pp-wave metrics correspond to
\begin{equation}\label{eq: H3 gR pp-waves}
[e_1,e_2] = \lambda_3 e_3,
\quad
[e_1,e_4] = \gamma_3 e_3,
\quad
[e_2, e_4] = \gamma_5 e_3.
\end{equation}
Now   the change of basis
$\bar e_1 = -\tfrac{\gamma_5}{\lambda_3} e_1+\tfrac{\gamma_3}{\lambda_3}e_2-e_4$,\,\,
$\bar e_2= \lambda_3 e_3$,\,\,
$\bar e_3 = -  e_2$,\,\,
$\bar e_4 = e_1$,\,\,
transforms the Lie bracket into
$[\bar e_3, \bar e_4] = \bar e_2$.
Therefore we conclude that left-invariant metrics \eqref{eq: H3 gR pp-waves} are also realized on $\algsdR$.

 	\begin{remark}\rm
 	Any left-invariant metric \eqref{H3-riemannian} with harmonic curvature is also realized on $\algsdR$. Indeed, working in the same polynomial ring as above, the computation of a Gröbner basis for the ideal 
 		$\mathcal{I}
 		= \langle \Cotton_{ijk} \cup \{\lambda_3\lambda'_3-1 \} \rangle$ leads to $9$ polynomials, three of them being $\mathbf{g}_1=\gamma_2(4\gamma_4^2+\lambda_3^2)$,
 		$\mathbf{g}_2=\gamma_4(4\gamma_4^2+\lambda_3^2)$ and $\mathbf{g}_3=\gamma_1-\gamma_4$, which imply $\gamma_1=\gamma_2=\gamma_4=0$ as above.
 	\end{remark}

\subsection{Semi-direct extensions with Lorentzian normal subgroup $\bm{\mathcal{H}^3$}}\label{H3-se:EGR-Lorentzian}
  It was shown by Rahmani \cite{Rahmani} that there are three non-homothetic classes of left-invariant Lorentz metrics on $\mathcal{H}^3$ corresponding to structure operator of types Ia and II. If the structure operator $L$ is diagonalizable, one distinguishes the two cases corresponding to $\operatorname{ker}L$ being positive definite or of Lorentzian signature. If the structure operator is of type II, then it is necessarily nilpotent.
We shall now analyze   these three cases separately.

\subsubsection{\bf The structure operator is diagonalizable of rank one with positive definite kernel}\label{H3-se:Lorentz-1}

In this case there exists an orthonormal basis $\{e_1,e_2,e_3,e_4\}$ of $\mathfrak{g}=\algsdH$, with $e_3$ timelike, so that 
\begin{equation}\label{H3-Ia+}
\begin{array}{ll}
[e_1,e_2]=-\lambda_3 e_3,  &
[e_1,e_4]=\gamma_1 e_1 - \gamma_2 e_2 + \gamma_3 e_3,  
\\
\noalign{\medskip}
[e_2,e_4]=\gamma_2 e_1 + \gamma_4  e_2 + \gamma_5 e_3,  &
[e_3,e_4]=(\gamma_1 +\gamma_4)e_3,
\end{array} 
\end{equation}
where $\lambda_3\neq 0$ and $\gamma_i\in\mathbb{R}$.

 \medskip

 We take an auxiliary variable $\lambda'_3$ to indicate that $\lambda_3\neq 0$	by means of  the polynomial $\lambda_3 \lambda'_3-1$ and
   compute a Gröbner basis for the ideal $\mathcal{I}
 = \langle \RicciRiccicomp{i}{j} \cup \{\tau, \lambda_3\lambda'_3-1 \} \rangle$ in the  polynomial ring $\mathbb{R}[\lambda'_3,\gamma_1,\gamma_2,\gamma_3,\gamma_5,\gamma_4,\lambda_3]$. As a consequence  we  obtain $20$ polynomials, two of them being
$
 \mathbf{g}_1 = (4\gamma_4^2+\lambda_3^2)(16\gamma_4^2+\lambda_3^2)\gamma_4^2$
and 
$ \mathbf{g}_2 =
			 \gamma_3^2+\gamma_5^2+\lambda_3^2+10\gamma_1\gamma_4$.
 Hence $\gamma_4=0$, which leads to  $\mathbf{g}_2\neq 0$ since  $\lambda_3\neq 0$. Therefore there is no left-invariant pp-wave metric in this case.

	\begin{remark}\rm
Observe that no left-invariant metric \eqref{H3-Ia+} has harmonic curvature. Indeed, a direct calculation gives 
		$\Cotton_{123}=\left( (\gamma_1+\gamma_4)^2 + \gamma_3^2 + \gamma_5^2 +\lambda_3^2\right)\lambda_3$, which is non-zero.
	\end{remark}

\subsubsection{\bf The structure operator is diagonalizable of rank one and with Lorentzian kernel}\label{H3-se:Lorentz-2}

In this setting, it is possible to choose    an orthonormal basis $\{e_1,e_2,e_3,e_4\}$ of $\mathfrak{g}=\algsdH$, with $e_3$ timelike, so that 
the left-invariant metrics are described by
\begin{equation}\label{H3-Ia-}
\begin{array}{ll}
[e_1,e_3]=-\lambda_2 e_2,  &
[e_1,e_4]=\gamma_1 e_1+\gamma_2 e_2 + \gamma_3 e_3,  
\\
\noalign{\medskip}
[e_2,e_4]=\gamma_4 e_2,  &
[e_3,e_4]=\gamma_5 e_1+\gamma_6 e_2-(\gamma_1 -\gamma_4)e_3,
\end{array}  
\end{equation}
where $\lambda_2\neq 0$ and $\gamma_i\in\mathbb{R}$.

\medskip

We proceed as in the previous cases. Let $\lambda_2'$ be an auxiliary variable so that the polynomial $\lambda_2 \lambda_2'-1$ reflects the non-vanishing of $\lambda_2$. We consider the ideal $\mathcal{I}= \langle \RicciRiccicomp{i}{j} \cup \{\tau, \lambda_2\lambda'_2-1 \} \rangle$ in the polynomial ring $\mathbb{R}[\lambda_2',\gamma_2,\gamma_6,\gamma_1,\gamma_3,\gamma_4,\gamma_5,\lambda_2]$. Computing a Gröbner basis we get a set of $48$ polynomials containing
\[
\begin{array}{l}
	\mathbf{g_1}= (16\gamma_4^2+\lambda_2^2)\gamma_4^3,\quad
	\mathbf{g_2}= (16\gamma_4^2+\lambda_2^2)\gamma_5^3,\quad 
	\mathbf{g_3}= (\gamma_3+\gamma_5)(\gamma_3^2+\gamma_5^2-\gamma_3\gamma_5),
	\\ \noalign{\medskip}
	\mathbf{g_4}=20 \gamma_1^2 - 7 \gamma_3^2 - 40 \gamma_4^2 - 7 \gamma_5^2 - 
	20 \gamma_1 \gamma_4 + 6 \gamma_3 \gamma_5 .
\end{array}
\]
Since $\lambda_2\neq 0$ it follows that  $\gamma_1=\gamma_3=\gamma_4=\gamma_5=0$ 
and the possible left-invariant pp-wave metrics are given by
\begin{equation}\label{eq: H3 gL.Ia- pp-waves}
	[e_1,e_3] = -\lambda_2 e_2,
	\quad
	[e_1,e_4] = \gamma_2 e_2,
	\quad
	[e_3, e_4] = \gamma_6 e_2.
\end{equation}
Note that   the change of basis
$\bar e_1 = -\tfrac{\gamma_6}{\lambda_2} e_1+\tfrac{\gamma_2}{\lambda_2}e_3+e_4$,\,\,
$\bar e_2= \lambda_2 e_2$,\,\,
$\bar e_3 =   e_1$,\,\,
$\bar e_4 = -e_3$,
%
transforms the Lie bracket into $[\bar e_3, \bar e_4] = \bar e_2$.
Therefore we conclude that left-invariant metrics given by Equation~\eqref{eq: H3 gL.Ia- pp-waves} are also realized on $\algsdR$.

	\begin{remark}\rm
Any Cotton-flat left-invariant metric \eqref{H3-Ia-} is  also realized on $\algsdR$. Indeed, working in the same polynomial ring as above we compute a Gröbner basis for the ideal 
		$\mathcal{I}= \langle \Cotton_{ijk} \cup \{\lambda_2\lambda'_2-1 \} \rangle$ and get $19$ polynomials, among which we find
		$\mathbf{g}_1=\gamma_3(\gamma_4^2+\lambda_2^2)$,
		$\mathbf{g}_2=\gamma_4(\gamma_4^2+\lambda_2^2)$,
		$\mathbf{g}_3=\gamma_5(\gamma_4^2+\lambda_2^2)$
		and $\mathbf{g}_4=4\gamma_1^2-\gamma_4^2-2\gamma_3(\gamma_3-\gamma_5)$. Hence,  $\gamma_1=\gamma_3=\gamma_4=\gamma_5=0$ as above. 
	\end{remark}

\subsubsection{\bf  The structure operator is 2-step nilpotent}\label{H3-se:Lorentz-3} 

In this   case  there exists a pseudo-orthonormal basis $\{u_1,u_2,u_3,u_4\}$ of $\mathfrak{g}=\algsdH$, with $\langle u_1,u_2\rangle=\langle u_3,u_3\rangle=\langle u_4,u_4\rangle=1$,  so that 
\begin{equation}\label{H3-II}
\begin{array}{ll}
[u_1,u_3]=-\varepsilon u_2,&
[u_1,u_4]=\gamma_1 u_1+\gamma_2 u_2+\gamma_3 u_3,   

\\
\noalign{\medskip}

[u_2,u_4]=\gamma_4 u_2, &
[u_3,u_4]=\gamma_5 u_1+\gamma_6 u_2 - (\gamma_1-\gamma_4)u_3 ,
\end{array} 
\end{equation}
with $\varepsilon^2=1$ and $\gamma_i\in\mathbb{R}$.

\medskip 

We consider the polynomial ring $\mathbb{R}[\ve, \gamma_2,\gamma_3,\gamma_6,\gamma_1,\gamma_4,\gamma_5]$. For the ideal  $\mathcal{I}$ generated by $\RicciRiccicomp{i}{j} \cup \{\tau, \ve^2-1 \}$ the computation of a Gröbner basis provides $52$ polynomials, among which we find
$\mathbf{g_1}=  \gamma_1^5$,
$\mathbf{g_2}=  \gamma_4^4$ and
$\mathbf{g_3}=  \gamma_5^3$.
Hence $\gamma_1=\gamma_4=\gamma_5=0$ so that
the possible left-invariant pp-wave metrics are given by
\begin{equation}\label{eq: H3 gL.II pp-waves}
	[u_1,u_3] = -\ve u_2,
	\quad
	[u_1,u_4] = \gamma_2 u_2+\gamma_3 u_3,
	\quad
	[u_3, u_4] = \gamma_6 u_2.
\end{equation}
If $\gamma_3\neq 0$ then the change of basis
$\bar u_1 = -\ve \gamma_6 u_1+\ve \gamma_2 u_3+u_4$,\,\,
$\bar u_2= -\ve  u_2$,\,\,
$\bar u_3 = \gamma_3   u_3$,\,\,
$\bar u_4 = -u_1$,
transforms the Lie bracket into
$[\bar u_1, \bar u_4] = \bar u_3$,\,\,
$[\bar u_3,\bar u_4] = \gamma_3 \bar u_2$,
while if $\gamma_3=0$ we consider
$\bar u_1 = -\ve \gamma_6 u_1+\ve \gamma_2 u_3+u_4$,\,\,
$\bar u_2= \ve  u_2$,\,\,
$\bar u_3 =    u_1$,\,\,
$\bar u_4 = -u_3$,
and the   Lie bracket transforms into $[\bar u_3,\bar u_4] = \bar u_2$.
Therefore   left-invariant metrics given by Equation~\eqref{eq: H3 gL.II pp-waves} are also realized on $\algsdR$.

	\begin{remark}\rm
Next we show that left-invariant metrics \eqref{H3-II} with harmonic curvature which are not realized on $\algsdR$ are necessarily locally conformally flat. Since $\Cotton_{123}=\frac{3\ve}{4}\gamma_5^2$ and 
		$\Cotton_{131}=-\frac{\ve}{4}\left( 4(3\gamma_1-\gamma_4)\gamma_4-\gamma_5(\gamma_3+7\gamma_6) \right)$, it follows that $\gamma_5=0$ and either $\gamma_4=0$ or $\gamma_4=3\gamma_1\neq 0$.
		If $\gamma_5=\gamma_4=0$ then $\Cotton_{142} = -\frac{3}{4}\gamma_1^3$, so that $\gamma_1=0$ and the  metric is also realized on $\algsdR$. Finally, if $\gamma_5=0$ and $\gamma_4=3\gamma_1\neq 0$, then $\Cotton_{143}=\frac{5}{2}\gamma_1^2(\gamma_3+\gamma_6)$ implies $\gamma_6=-\gamma_3$ and the space is locally conformally flat. 
	\end{remark}

\subsection{Semi-direct extensions with degenerate normal subgroup $\bm{\mathcal{H}^3}$}\label{H3-se:2}

Let $\mathfrak{g}=\algsdH$ be a four-dimensional Lie algebra with a Lorentzian inner product $\ip$  which restricts to a degenerate inner product on   $\mathfrak{h}_3$. We  consider the two-distinct situations when the restriction of the metric to the derived subalgebra $\mathfrak{h}'_3$ is degenerate or spacelike.

\subsubsection{\bf $\bm{{\mathfrak{h}_3'=\operatorname{\bf span}\{v\}}}$ is a null subspace}\label{H3-se:deg-null subspace}

In this case  setting $u_3=v$ we can take a pseudo-orthonormal basis $\{u_1,u_2,u_3,u_4\}$ of $\mathfrak{g}=\algsdH$,  with $\langle u_1,u_1\rangle=\langle u_2,u_2\rangle=\langle u_3,u_4\rangle=1$,  so that  
\begin{equation}\label{H3-eq: g_D0 ip} 
\begin{array}{ll}
[u_1,u_2]= \lambda_1 u_3 , &
[u_1,u_4]= \gamma_1 u_1 -  \gamma_2 u_2 + \gamma_3 u_3 ,
\\
\noalign{\medskip}
[u_2,u_4]=  \gamma_2 u_1 + \gamma_4 u_2 + \gamma_5 u_3 , &
[u_3,u_4]=(\gamma_1+\gamma_4) u_3 ,
\end{array}
\end{equation}
where  $\lambda_1\neq 0$ and    $\gamma_i\in\mathbb{R}$.

\medskip

A direct calculation shows that, with no additional assumptions, the only non-zero component of the Ricci tensor is given by $\rho(u_4,u_4)= \frac{1}{2}(4\gamma_1\gamma_4+\lambda_1^2)$, which proves that the Ricci operator is isotropic.
Moreover,   $\nabla_{u_i} u_3=0$ for $i=1,2,3$ and $\nabla_{u_4} u_3 = -(\gamma_1+\gamma_4)u_3$, so that $u_3$ determines   a null recurrent  left-invariant vector field. Furthermore, a straightforward calculation shows that $R(x,y)=0$ and $\nabla_xR=0$ for all $x,y\in u_3^\perp=\span\{u_3,u_1,u_2\}$. 
Hence left-invariant metrics determined by \eqref{H3-eq: g_D0 ip} are plane waves (see~\cite{Leistner}). Finally, since $\lambda_1\neq 0$, we may take $\lambda_1=1$ working in the homothetic class of the initial metric,  corresponding to the metrics given in Theorem~\ref{H3-th: pp-waves}.

\subsubsection{\bf $\bm{{\mathfrak{h}_3'=\span\{v\}}}$ is a spacelike subspace}\label{H3-se:deg-spacelike subspace}

We set  $u_1=\frac{v}{\|v\|}$ and consider a pseudo-orthonormal basis $\{u_1,u_2,u_3,u_4\}$ of $\mathfrak{g}=\algsdH$, with $\langle u_1,u_1\rangle=\langle u_2,u_2\rangle=\langle u_3,u_4\rangle=1$,  so that  
\begin{equation}\label{H3-degenerate-spacelike}
\begin{array}{l}
[u_1,u_4]= \gamma_1 u_1, \quad
[u_2,u_3]= \lambda_3 u_1, \quad
[u_2,u_4]= \gamma_2 u_1 +\gamma_3 u_2 + \gamma_4 u_3,
\\
\noalign{\medskip}
[u_3,u_4]=  \gamma_5 u_1+\gamma_6 u_2 + (\gamma_1-\gamma_3) u_3 ,
\end{array} 
\end{equation}
where  $\lambda_3\neq 0$ and $\gamma_i\in\mathbb{R}$.

\medskip

We start considering an   auxiliary variable $\lambda'_3$ to indicate that $\lambda_3\neq 0$	by means of  $\lambda_3 \lambda'_3-1$, and   the ideal $\mathcal{I}= \langle \RicciRiccicomp{i}{j} \cup \{\tau, \lambda_3\lambda'_3-1 \} \rangle$ in the polynomial ring $\mathbb{R}[\lambda_3',\lambda_3,\gamma_2,\gamma_5,\gamma_6,\gamma_4,\gamma_3,\gamma_1]$. Computing a Gröbner basis for this ideal we get $27$ polynomials, four of them being
$\mathbf{g_1}=  \gamma_1^3$,
$\mathbf{g_2}=  2\gamma_3^2-3\gamma_1^2+2\gamma_1\gamma_3$,
$\mathbf{g_3}=  \gamma_4^3$ and
$\mathbf{g_4}=  \gamma_6^2$.
Hence $\gamma_1=\gamma_3=\gamma_4=\gamma_6=0$ 
and the possible left-invariant pp-wave metrics are given by
\begin{equation}\label{eq: H3 gD+ pp-waves}
	[u_2,u_3] = \lambda_3 u_1,
	\quad
	[u_2,u_4] = \gamma_2 u_1,
	\quad
	[u_3, u_4] = \gamma_5 u_1.
\end{equation}
Now    the change of basis
$\bar u_1 = \tfrac{\gamma_5}{\lambda_3} u_2-\tfrac{\gamma_2}{\lambda_3}u_3+u_4$,\,\,
$\bar u_2= \lambda_3 u_1$,\,\,
$\bar u_3 =   u_3$,\,\,
$\bar u_4 = -u_2$,
transforms the Lie bracket into $[\bar u_3, \bar u_4] = \bar u_2$.
Therefore we conclude that left-invariant metrics given by Equation~\eqref{eq: H3 gD+ pp-waves} are also realized on $\algsdR$.

	\begin{remark}\rm
Any Cotton-flat left-invariant metric \eqref{H3-degenerate-spacelike} is  also realized on $\algsdR$. Indeed, considering  the same polynomial ring as above and computing a Gröbner basis for  
		$\mathcal{I}= \langle \Cotton_{ijk} \cup \{\lambda_3\lambda'_3-1 \} \rangle$, one obtains   $7$ polynomials including
		$\mathbf{g}_1=\gamma_1$,
		$\mathbf{g}_2=\gamma_3$,
		$\mathbf{g}_3=\gamma_4$
		and $\mathbf{g}_4=\gamma_6$. Hence, $\gamma_1=\gamma_3=\gamma_4=\gamma_6=0$ as above.
	\end{remark}

\section{Semi-direct extensions of the Euclidean and Poincaré Lie groups}\label{se:EE}

If a semi-direct extension of the Euclidean or Poincaré Lie groups is unimodular then 
the results in  \cite{ABDO} imply that it reduces to a direct product $E(1,1)\times \mathbb{R}$ or $\widetilde{E}(2)\times\mathbb{R}$, and thus it is  is isomorphic to a semi-direct extension $\mathbb{R}^3\rtimes\mathbb{R}$ of the Abelian Lie group covered by the analysis in Section~\ref{se:R3}. Taking this fact into account,   the use of Gröbner bases becomes in this case a fundamental tool to prove the following result, which is obtained as a consequence of the subsequent analysis.

\begin{theorem}\label{E11-E2-th: pp-waves}
	A four-dimensional  semi-direct extension of the Euclidean or the Poincaré Lie groups which is a non-unimodular pp-wave is 
	isomorphically homothetic to a left-invariant plane wave metric on $\operatorname{Aff}(\mathbb{C})$ determined by
	\[
	\begin{array}{lll}
		[v_1,v_3]=  v_2,  &
		[v_2,v_3]=-  v_1,
		\\ \noalign{\medskip}
		[v_1,v_4]=\gamma_1 v_1 + \gamma_2 v_2, & 
		[v_2,v_4]=-\gamma_2 v_1+\gamma_1 v_2,
		\qquad\gamma_1\neq 0,
	\end{array}
	\]
	
	\smallskip
	
	\noindent
where $\{v_i\}$ is a pseudo-orthonormal basis   with $\langle v_1,v_1\rangle=\langle v_2,v_2\rangle=\langle v_3,v_4\rangle=1$.

Moreover non-unimodular semi-direct extensions of the Euclidean or Poincaré Lie groups have harmonic curvature if and only if they are locally symmetric.
\end{theorem}

\begin{remark}\rm
Non-unimodular extensions of the Euclidean and Poincaré Lie groups are isomorphic to the affine groups determined by the Lie algebras $\mathfrak{aff}(\mathbb{C})$ and $\mathfrak{aff}(\mathbb{R})\times\mathfrak{aff}(\mathbb{R})$, respectively. Hence it also follows from the previous theorem that there are no left-invariant pp-wave metrics on $\mathfrak{aff}(\mathbb{R})\times\mathfrak{aff}(\mathbb{R})$.
	
		Left-invariant metrics given in Theorem~\ref{E11-E2-th: pp-waves} are locally conformally flat and locally symmetric but never Einstein. Hence they are plane waves of type~(i).
	
\end{remark}

\subsection{Extensions of the Riemannian Euclidean and Poincaré Lie groups}
Let $G=G_3\rtimes\mathbb{R}$ with $G_3=E(1,1)$ or $G_3=\widetilde{E}(2)$ so that the restriction of the metric to $G_3$ is positive definite. In this case  there exists an orthonormal basis $\{ e_1,e_2,e_3,e_4\}$ of the Lie algebra $\mathfrak{g}=\mathfrak{g}_3\rtimes\mathbb{R}$, with $e_4$ timelike, so that 
\begin{equation}\label{E-Riemann}
\begin{array}{lll}
[e_1,e_3]\!=\!-\lambda_2 e_2, & [e_2,e_3]\!=\!\lambda_1 e_1,& 
\\ \noalign{\medskip}
{[e_1,e_4]}\!=\!\gamma_1 e_1+\gamma_2\lambda_2 e_2, &{[e_2,e_4]}\!=\!-\gamma_2\lambda_1 e_1+\gamma_1 e_2,&[e_3,e_4]\!=\!\gamma_3 e_1+\gamma_4 e_2,
\end{array} 
\end{equation}
where $\lambda_1 \lambda_2\neq 0$ and $\gamma_i\in\mathbb{R}$.

A direct calculation shows that  the Lie algebra is unimodular if and only if $\gamma_1=0$. 
In the non-unimodular case  we introduce an auxiliary variable $\gamma'_1$ and use the polynomial $\gamma_1\gamma'_1-1$ to indicate that   $\gamma_1\neq 0$. Considering the polynomial ring 
$\mathbb{R}[\gamma_2,\gamma_3,\gamma_4,\lambda_1,\lambda_2,\gamma'_1,\gamma_1]$  and the ideal
$\mathcal{I}= \langle \RicciRiccicomp{i}{j} \cup \{\tau, \gamma_1\gamma'_1 - 1 \} \rangle$ 
we compute a Gröbner basis  and obtain $9$ polynomials, one of them being
\[
\mathbf{g} = (\gamma_1^2 + \lambda_2^2) (9 \gamma_1^2 + \lambda_2^2) 
		(64 \gamma_1^4 + 48 \gamma_1^2 \lambda_2^2 + \lambda_2^4).
\]
Clearly $\mathbf{g}\neq 0$, so that  we conclude the non-existence of  left-invariant  pp-wave metrics in the non-unimodular setting.

	\begin{remark}\label{re: E11-E2 Cotton-flat}\rm
	In order to analyze the Cotton-flat condition, we restrict to the non-unimodular case and consider the same polynomial ring as above, computing a Gröbner basis for the ideal 
		$\mathcal{I}= \langle \Cotton_{ijk} \cup \{\gamma_1\gamma'_1 - 1 \} \rangle$. One gets $5$ polynomials, among which we have
		$\mathbf{g}_1=\gamma_3$, 
		$\mathbf{g}_2=\gamma_4$, 
		$\mathbf{g}_3=\lambda_1^2-\lambda_2^2$ and
		$\mathbf{g}_4=-(\lambda_1-\lambda_2)(\gamma_1^2-(\gamma_2^2-1)\lambda_2^2)$.
		Hence, $\gamma_3=\gamma_4=0$ and either $\lambda_1=\lambda_2$ or
		$\lambda_1=-\lambda_2$, $\gamma_1^2-(\gamma_2^2-1)\lambda_2^2=0$. Finally, a direct calculation shows that in any of these cases the left-invariant metric \eqref{E-Riemann} is locally symmetric. 
	\end{remark}

\subsection{Extensions of the Lorentzian Euclidean and Poincaré Lie groups}

We consider the different possibilities for the Jordan normal form of the structure operator $L$.

\subsubsection{\bf Diagonalizable structure operator with spacelike $\bm{\ker L}$}\label{E-Ia-spacelike}

In this case  the left-invariant metrics are described in an orthonormal basis $\{ e_1,e_2,e_3,e_4\}$, with $e_3$ timelike, by
\begin{equation}\label{E-Lorentz-Ia+}
\begin{array}{lll}
[e_1,e_2]\!=\!-\lambda_3 e_3, & [e_1,e_3]\!=\!-\lambda_2 e_2, &
\\ \noalign{\medskip}
{[e_1,e_4]}\!=\!\gamma_1 e_2+\gamma_2 e_3, & [e_2,e_4]\!=\!\gamma_3 e_2+\gamma_4\lambda_3 e_3,&[e_3,e_4]\!=\!\gamma_4 \lambda_2 e_2+\gamma_3 e_3,
\end{array}  
\end{equation}
where $\lambda_2\lambda_3\neq 0$ and $\gamma_i\in\mathbb{R}$.

The Lie algebra is unimodular if and only if $\gamma_3=0$. As in the previous case we consider the non-unimodular setting and introduce an auxiliary variable $\gamma'_3$  to indicate that   $\gamma_3\neq 0$ by means of   the polynomial $\gamma_3\gamma'_3 - 1$.
Let $\mathcal{I}$ be the ideal generated by $\RicciRiccicomp{i}{j} \cup \{\tau, \gamma_3\gamma'_3 - 1 \}$ in the polynomial   ring  
$\mathbb{R}[\gamma_1,\gamma_2,\gamma_4,\lambda_2,\lambda_3,\gamma'_3,\gamma_3]$.
The computation of a  Gröbner basis for this ideal  yields  $11$ polynomials, one of them being
$\mathbf{g} =  (\gamma_3^2 + \lambda_3^2) (9 \gamma_3^2 + \lambda_3^2) 
			(64 \gamma_3^4 + 48 \gamma_3^2 \lambda_3^2 + \lambda_3^4)$.
Since $\lambda_3\neq 0$ we have $\mathbf{g}\neq 0$ and  we conclude the non-existence of  left-invariant  pp-wave metrics  in the non-unimodular case. 

	\begin{remark}\rm
Any non-unimodular Cotton-flat left-invariant metric \eqref{E-Lorentz-Ia+} is locally symmetric.
		Indeed, considering the ideal $\mathcal{I}= \langle \Cotton_{ijk} \cup \{\gamma_3\gamma'_3 - 1 \} \rangle$ in the same polynomial ring as above, we compute a Gröbner basis to obtain a set of $5$ polynomials containing 
		$\mathbf{g}_1=\gamma_1$, 
		$\mathbf{g}_2=\gamma_2$ and
		$\mathbf{g}_3=(\lambda_2-\lambda_3)(\gamma_3^2+(\gamma_4^2+1)\lambda_3^2)$. Thus, $\gamma_1=\gamma_2=0$ and $\lambda_2=\lambda_3$, from where it follows that the metric is locally symmetric. 
	\end{remark}

\subsubsection{\bf Diagonalizable structure operator with timelike $\bm{\ker L}$}\label{E-Ia-timelike}

The corresponding left-invariant metrics are described in an orthonormal basis $\{ e_1,e_2,e_3,e_4\}$, with $e_3$ timelike, by
\begin{equation}\label{E-Lorentz-Ia-}
\begin{array}{lll}
[e_1,e_3]\!=\!-\lambda_2 e_2, & [e_2,e_3]\!=\!\lambda_1 e_1, &
\\ \noalign{\medskip}
{[e_1,e_4]}\!=\!\gamma_1 e_1+\gamma_2\lambda_2 e_2, & [e_2,e_4]\!=\!-\gamma_2\lambda_1 e_1+\gamma_1 e_2, &[e_3,e_4]\!=\!\gamma_3 e_1+\gamma_4 e_2,
\end{array} 
\end{equation}
where $\lambda_1\lambda_2\neq 0$ and   $\gamma_i\in\mathbb{R}$.

One easily checks that the Lie algebra is unimodular if and only if $\gamma_1=0$.
Assuming $\gamma_1\neq 0$ and adding an auxiliary variable $\gamma_1'$ to use the polynomial  $\gamma_1\gamma'_1 - 1$, we consider the polynomial ring $\mathbb{R}[\gamma_2,\gamma_3,\gamma_4,\lambda_1,\lambda_2,\gamma'_1,\gamma_1]$. We compute a Gröbner basis for the ideal 
$\mathcal{I}= \langle \RicciRiccicomp{i}{j} \cup \{\tau, \gamma_1\gamma'_1 - 1 \} \rangle$ and get a set of $9$ polynomials, among which we find
$\mathbf{g}=(\gamma_1^2 + \lambda_2^2) (9 \gamma_1^2 + \lambda_2^2) 
		(64 \gamma_1^4 + 48 \gamma_1^2 \lambda_2^2 + \lambda_2^4)$.
Since $\lambda_2\neq 0$ it follows that $\mathbf{g}\neq 0$, and therefore  we conclude the non-existence of  left-invariant  pp-wave metrics  in   the non-unimodular setting.

	\begin{remark}\rm
		We proceed exactly as in Remark~\ref{re: E11-E2 Cotton-flat} to get that, in the non-unimodular case, any Cotton-flat left-invariant metric \eqref{E-Lorentz-Ia-} is locally symmetric. 
	\end{remark}

\subsubsection{\bf Structure operator $\bm{L}$ with a complex eigenvalue}\label{E-complex}

The left-invariant metrics are determined by 
\begin{equation}\label{E-Lorentz-Ib}
\begin{array}{l}
[e_1,e_2]\!=\!-\beta e_2-\alpha e_3, 
\qquad 
[e_1,e_3]\!=\!-\alpha e_2+\beta e_3,\qquad [e_1,e_4]\!=\!\gamma_1 e_2+\gamma_2 e_3,
\\ \noalign{\medskip}
[e_2,e_4]\!=\!2\gamma_3\beta e_2+(\gamma_3-\gamma_4)\alpha e_3,
\quad
[e_3,e_4]\!=\!(\gamma_3-\gamma_4)\alpha e_2+2\gamma_4\beta e_3 ,
\end{array} 
\end{equation}
where $\beta\neq 0$, $\alpha,\gamma_i\in\mathbb{R}$ and
$\{ e_1,e_2,e_3,e_4\}$ is an orthonormal basis with $e_3$ timelike. 

In this case the computation of a Gröbner basis for the ideal 
$\mathcal{I}= \langle \RicciRiccicomp{i}{j} \cup \{\tau\} \rangle$ in the polynomial ring 
$\mathbb{R}[\alpha,\gamma_2,\gamma_3,\gamma_4,\gamma_1,\beta]$ provides a set of $74$ elements containing the polynomial
$\mathbf{g} = 
	  \left( (8 \gamma_4^2 + 2) \beta^2 + \gamma_1^2 \right)^2 
	  \left( (\beta^2 + 2 \gamma_1^2) \beta^2 + 2 \gamma_1^4\right) 
	  (8 \beta^2 + \gamma_1^2) \beta$.
Note that $\mathbf{g}$ does not vanish since $\beta\neq 0$. Hence   we conclude that left-invariant metrics in this case are never pp-waves.

	\begin{remark}\rm
Any Cotton-flat left-invariant metric \eqref{E-Lorentz-Ib} is locally symmetric. Indeed, considering the component $\Cotton_{144}=-\frac{3}{4}(\gamma_1^2+\gamma_2^2)\beta$, one has that $\gamma_1=\gamma_2=0$.
		Now, $\Cotton_{231}=4 \left( (\gamma_3-\gamma_4)^2+1 \right)\alpha\beta^2$, and therefore $\alpha=0$. With these conditions a direct calculation shows that the vanishing of the Cotton tensor is determined by the component $\Cotton_{122}=-2(4\gamma_3\gamma_4-1)\beta^3$. Hence, $4\gamma_3\gamma_4=1$ and one easily checks that the corresponding metric is locally symmetric. 
	\end{remark}

\subsubsection{\bf Structure operator of type II with degenerate kernel}\label{E-II-degenerate}
In this case the left-invariant metrics are described, in terms of a pseudo-orthonormal basis $\{ u_1,u_2,u_3,u_4\}$ of the Lie algebra with $\langle u_1,u_2\rangle=\langle u_3,u_3\rangle=\langle u_4,u_4\rangle=1$, by the Lie brackets
\begin{equation}\label{E-Lorentz-II-degenerate}
\begin{array}{lll}
[u_1,u_2]\!=\!\lambda_2 u_3,& [u_1,u_3]\!=\!{-\varepsilon} u_2, &
\\ \noalign{\medskip}
[u_1,u_4]\!=\!\gamma_1 u_2+\gamma_2 u_3, & 
[u_2,u_4]\!=\!\gamma_3 u_2+\gamma_4\lambda_2  u_3, & 
[u_3,u_4]\!=\!\gamma_3 u_3-\varepsilon\gamma_4 u_2,
\end{array}
\end{equation}
where $\lambda_2\neq 0$, $\ve=\pm 1$ and  $\gamma_i\in\mathbb{R}$. 

Let $\mathcal{I}$ be the ideal generated by $\RicciRiccicomp{i}{j} \cup \{\tau, \ve^2-1\} $ in the polynomial ring 
$\mathbb{R}[\gamma_1,\gamma_2,\gamma_3,\gamma_4,\ve,\lambda_2]$. We compute a Gröbner basis for this ideal and get $38$ polynomials, one of them being 
$\mathbf{g} = \lambda_2^5\neq 0$.
Hence, as in the previous case, left-invariant metrics in this case are never pp-waves.

	\begin{remark}\rm
		To analyze the vanishing of the Cotton tensor   we calculate the components $\Cotton_{344}=-\frac{1}{4}\ve\gamma_4^2\lambda_2^2$ and 
		$\Cotton_{132}=\frac{1}{2} (\gamma_3^2+\lambda_2^2-2\gamma_2\gamma_4\lambda_2)  \lambda_2$. Hence, 
		$\gamma_4=0$ and $\Cotton_{132}=\frac{1}{2} (\gamma_3^2+\lambda_2^2)  \lambda_2$ does not vanish since $\lambda_2\neq 0$, thus showing that no left-invariant metric \eqref{E-Lorentz-II-degenerate} may have harmonic curvature.
	\end{remark}

\subsubsection{\bf Structure operator of type II with spacelike kernel}\label{E-II-spacelike}
The left-invariant metrics are described by
\begin{equation}\label{E-Lorentz-II+}
\begin{array}{lll}
	[u_1,u_3]\!=\! -\lambda_1 u_1 \!-\!\varepsilon u_2 , \!&\!
	[u_2,u_3]\!=\!\lambda_1 u_2 ,
	\\ \noalign{\medskip}
	[u_1,u_4]\!=\!\gamma_1 u_1\!+\!\gamma_2 u_2, \!&\!
	[u_2,u_4]\!=\!(\gamma_1\!-\!2\varepsilon\gamma_2\lambda_1)u_2, \!& \!
	[u_3,u_4]\!=\!\gamma_3 u_1\!+\!\gamma_4 u_2,
\end{array}
\end{equation}
where $\lambda_1\neq 0$, $\ve=\pm 1$, $\gamma_i\in\mathbb{R}$ and 
$\{ u_1,u_2,u_3,u_4\}$ is a pseudo-orthonormal basis of the Lie algebra with $\langle u_1,u_2\rangle=\langle u_3,u_3\rangle=\langle u_4,u_4\rangle=1$. 

A direct calculation shows that the Lie algebra is unimodular if and only if $\varepsilon\gamma_2\lambda_1-\gamma_1=0$. We consider the non-unimodular case and introduce an auxiliary variable $\varphi$ to indicate that $\ve\gamma_2\lambda_1-\gamma_1\neq 0$ using the polynomial $(\ve \gamma_2\lambda_1-\gamma_1)\varphi - 1$.
In the polynomial ring 
$\mathbb{R}[\varphi, \ve,\lambda_1,\gamma_1,\gamma_2,\gamma_3,\gamma_4]$ we compute a Gröbner basis for the ideal $\mathcal{I}=
\langle \RicciRiccicomp{i}{j} \cup \{\tau, (\varepsilon\gamma_2\lambda_1-\gamma_1)\varphi - 1\}  \rangle$ and   observe that it reduces to $\{1\}$. 
Thus   we conclude the non-existence of  left-invariant  pp-wave metrics  in the non-unimodular case.

	\begin{remark}\rm
Left-invariant metrics \eqref{E-Lorentz-II+} are never Cotton-flat. Indeed, vanishing of
		$\Cotton_{232}=-\frac{5}{4}\gamma_3^2\lambda_1$ implies $\gamma_3=0$, and then
		$\Cotton_{134} = \frac{1}{2}\gamma_4(\gamma_1^2-\ve\gamma_1\gamma_2\lambda_1-\lambda_1^2)$, while
		$\Cotton_{341} = \frac{1}{2}\gamma_4(2\gamma_1^2-\ve\gamma_1\gamma_2\lambda_1+\lambda_1^2)$. 
		Hence, $\Cotton_{134}-\Cotton_{341}=-\frac{1}{2}\gamma_4(\gamma_1^2+2\lambda_1^2)$, so that $\gamma_4=0$. Finally, using $\gamma_3=\gamma_4=0$ we calculate 
		$\Cotton_{131}=-2\ve \left(\gamma_1^2+(\gamma_2^2+2)\lambda_1^2\right)$, which does not vanish.
		
	\end{remark}

\subsubsection{\bf Structure operator of type III}\label{E-III}
In this case the left-invariant metrics are described by
\begin{equation}\label{E-Lorentz-III}
\begin{array}{lll}
	[u_1,u_2]=u_1,  &
	[u_2,u_3]=u_3,
	\\ \noalign{\medskip}
	[u_1,u_4]=\gamma_1 u_1, & 
	[u_2,u_4]=\gamma_2 u_1+\gamma_3 u_3, & 
	[u_3,u_4]=\gamma_4 u_3,
\end{array}
\end{equation}
where $\gamma_i\in\mathbb{R}$ and  $\{ u_1,u_2,u_3,u_4\}$ is a pseudo-orthonormal basis of the Lie algebra with $\langle u_1,u_2\rangle=\langle u_3,u_3\rangle=\langle u_4,u_4\rangle=1$. 

One easily checks that the Lie algebra is unimodular if and only if $\gamma_1+\gamma_4=0$. 
In the non-unimodular case,  let $\varphi$ be an auxiliary variable which reflects the fact that $\gamma_1+\gamma_4\neq 0$ by means of the polynomial $(\gamma_1+\gamma_4)\varphi-1$. We work in the polynomial ring 
$\mathbb{R}[\varphi,\gamma_1,\gamma_2,\gamma_3,\gamma_4]$.
As in the previous case, the computation of a Gröbner basis for the ideal 
$\mathcal{I}$ generated by 
$\RicciRiccicomp{i}{j} \cup \{\tau,  (\gamma_1+\gamma_4)\varphi-1\}$  reduces to $\{1\}$. 
Therefore, the is no left-invariant pp-wave metric in the non-unimodular setting.

	\begin{remark}\rm
		In the non-unimodular case any Cotton-flat left-invariant metric \eqref{E-Lorentz-III} is locally symmetric.
		To show it, we work in the same polynomial ring as above and compute a Gröbner basis for the ideal 
		$\mathcal{I}= \langle \Cotton_{ijk} \cup \{ (\gamma_1+\gamma_4)\varphi-1 \} \rangle$. Thus, we get  $4$ polynomials, three of them being 
		$\mathbf{g}_1=\gamma_3$, 
		$\mathbf{g}_2=(\gamma_2\gamma_4+1)(\gamma_2\gamma_4+2)$ and
		$\mathbf{g}_3=\gamma_1+2(\gamma_2\gamma_4+1)\gamma_4$. 
		Hence, $\gamma_3=0$ and either $\gamma_1=0$, $\gamma_2\gamma_4=-1$ or
		$\gamma_1=2\gamma_4$ and $\gamma_2\gamma_4=-2$. Finally, a direct calculation shows that, in any case, the underlying metric is locally symmetric.  
	\end{remark}

\subsection{Extensions of degenerate Euclidean and Poincaré Lie groups}
Let $G=G_3\rtimes\mathbb{R}$ with $G_3=E(1,1)$ or $G_3=\widetilde{E}(2)$.
Since the restriction of the metric to $\mathfrak{g}_3$ is degenerate of signature $(++0)$, we consider separately the cases when the induced metric on 
the derived algebra $\mathfrak{g}_3'=[\mathfrak{g}_3,\mathfrak{g}_3]$ is Riemannian in Section~\ref{EE-se:riemannian}, while the case when the restriction of the metric to  $\mathfrak{g}'_3$ is degenerate in considered in Section~\ref{EE-se:degenerate}.

\subsubsection{\bf The induced metric on the derived algebra $\bm{\mathfrak{g}'_3}$ is positive definite}\label{EE-se:riemannian}

There exists a pseudo-orthonormal basis $\{u_1,u_2,u_3,u_4\}$ of the Lie algebra $\mathfrak{g}=\mathfrak{g}_3\rtimes\mathbb{R}$ with $\langle u_1,u_1\rangle=\langle u_2,u_2\rangle=\langle u_3,u_4\rangle=1$ so that $\mathfrak{g}_3=\operatorname{span}\{ u_1,u_2,u_3\}$ and the derived algebra $\mathfrak{g}'_3=\operatorname{span}\{ u_1,u_2\}$.
In this situation, using that $\mathfrak{g}_3$ is unimodular and the derived algebra $\mathfrak{g}'_3$ is Abelian, the Jacobi identity leads to two-distinct situations which we consider as follows.

\subsubsection*{\underline{Case 1}}
There exists a pseudo-orthonormal basis $\{u_1,u_2,u_3,u_4\}$ of the Lie algebra with $\langle u_1,u_1\rangle=\langle u_2,u_2\rangle=\langle u_3,u_4\rangle=1$ so that
\begin{equation}\label{eq:E-degenerate-+1}
\begin{array}{lll}
[u_1,u_3]=\lambda_1 u_2,  &\!
[u_2,u_3]=-\lambda_1 u_1,
\\ \noalign{\medskip}
[u_1,u_4]=\gamma_1 u_1 + \gamma_2 u_2, & \!
[u_2,u_4]=-\gamma_2 u_1+\gamma_1 u_2, & \!
[u_3,u_4]=\gamma_3 u_1 + \gamma_4 u_2,
\end{array}
\end{equation}
where $\lambda_1\neq 0$ and $\gamma_i\in\mathbb{R}$. 

A direct calculation shows that the Lie algebra $\mathfrak{g}=\mathfrak{g}_3\rtimes\mathbb{R}$ is unimodular if and only if $\gamma_1 =0$.  Next we analyze the non-unimodular case, i.e., $\gamma_1 \neq 0$.
A straightforward calculation shows that the Ricci operator  satisfies
$$
\begin{array}{lll}
\Riccicomp{1}{1} =  -\frac{1}{2} \gamma_3^2,
&
 \Riccicomp{3}{1}=\Riccicomp{1}{4}=
\frac{1}{2}(3\gamma_1\gamma_3+\gamma_2\gamma_4)  ,
&
\Riccicomp{4}{1}=\Riccicomp{1}{3}=
-\frac{1}{2}\gamma_4\lambda_1, 
\\ 
\noalign{\medskip}
\Riccicomp{2}{2} =   -\frac{1}{2} \gamma_4^2,
&
\Riccicomp{3}{2}=\Riccicomp{2}{4}=\frac{1}{2}(3\gamma_1\gamma_4-\gamma_2\gamma_3),	
&
\Riccicomp{4}{2}=\Riccicomp{2}{3}=
\frac{1}{2}\gamma_3\lambda_1,
\\ 
\noalign{\medskip}
\Riccicomp{3}{4}=-2\gamma_1^2,
&
\Riccicomp{3}{3} =  \Riccicomp{4}{4} = 
\frac{1}{2} (\gamma_3^2+\gamma_4^2),
&
\Riccicomp{2}{1} = \Riccicomp{1}{2} = -\frac{1}{2}\gamma_3\gamma_4.
\end{array}
$$
%
%
%
%
%
%
%
%
%
%
%
%
%
%
Hence, one easily checks that $\Ricci^2=0$ if and only if $\gamma_3=\gamma_4=0$, in which case the corresponding left-invariant metric is given by
\begin{equation}\label{eq: aff(C) plane wave}
	\begin{array}{lll}
		[u_1,u_3]=\lambda_1 u_2,  &
		[u_2,u_3]=-\lambda_1 u_1,
		\\ \noalign{\medskip}
		[u_1,u_4]=\gamma_1 u_1 + \gamma_2 u_2, & 
		[u_2,u_4]=-\gamma_2 u_1+\gamma_1 u_2.
	\end{array}
\end{equation}
Now a straightforward calculation shows that  $u_3$ determines a left-invariant null parallel vector field. Moreover, $R(x,y)=0$ and $\nabla_x R=0$ for all $x,y\in u_3^\perp=\span\{u_3,u_1,u_2\}$. 
Therefore, left-invariant metrics given by Equation~\eqref{eq: aff(C) plane wave} are plane waves (see~\cite{Leistner}). Note that in the non-unimodular case ($\gamma_1\neq 0$) these metrics are realized on the Lie group corresponding to $\mathfrak{aff}(\mathbb{C})$. Moreover, since $\lambda_1\neq 0$, we may take $\lambda_1=1$ working in the homothetic class of the initial metric,  obtaining  the metrics  given in Theorem~\ref{E11-E2-th: pp-waves}. A straightforward calculation shows  that left-invariant metrics determined by \eqref{eq: aff(C) plane wave} are locally symmetric and locally conformally flat.

 	\begin{remark}\rm
 		In the non-unimodular case the only Cotton-flat left-invariant metrics \eqref{eq:E-degenerate-+1} correspond to the plane waves given by~\eqref{eq: aff(C) plane wave}, since   $\Cotton_{123}=\frac{1}{4}(\gamma_3^2+\gamma_4^2)\lambda_1=0$ implies $\gamma_3=\gamma_4=0$.
 	\end{remark}

\subsubsection*{\underline{Case 2}}
There exists a pseudo-orthonormal basis $\{u_1,u_2,u_3,u_4\}$ of the Lie algebra, with $\langle u_1,u_1\rangle=\langle u_2,u_2\rangle=\langle u_3,u_4\rangle=1$, so that
\begin{equation}\label{eq:E-degenerate-+2}
\begin{array}{l}
[u_1,u_3]\!=\!\lambda_1 u_1+\lambda_2 u_2,  \quad
[u_2,u_3]\!=\! -\lambda_2 u_1-\lambda_1 u_2, \quad
[u_3,u_4]\!=\!\gamma_3 u_1 + \gamma_4 u_2,

\\ \noalign{\medskip}
[u_1,u_4]\!=\! \gamma_1 u_1 + \frac{(\gamma_1-\gamma_2)\lambda_2}{2\lambda_1} u_2, \quad
[u_2,u_4]\!=\! - \frac{(\gamma_1-\gamma_2)\lambda_2}{2\lambda_1} u_1+\gamma_2 u_2,  
\end{array}
\end{equation}
where $\lambda_1\neq 0$, $\lambda_1^2-\lambda_2^2\neq 0$ and  $\gamma_i\in\mathbb{R}$. 

One easily checks that, in this case,  the Lie algebra   is unimodular if and only if $\gamma_1+\gamma_2=0$. In the non-unimodular case we introduce an auxiliary variable $\varphi$ and use the polynomial 
$(\gamma_1+\gamma_2)\varphi -1$. 
In the polynomial ring $\mathbb{R}[\varphi,\gamma_1,\gamma_2,\gamma_3,\gamma_4,\lambda_1,\lambda_2]$ we consider  
$\mathcal{I}=\langle  \lambda_1^2 \RicciRiccicomp{i}{j} \cup \{\tau, (\gamma_1+\gamma_2)\varphi -1\}\rangle$, where we have taken $\lambda_1^2 \RicciRiccicomp{i}{j}$ to remove the variables in the denominators. Computing a Gröbner basis for this ideal  we get 
a set of $91$ elements containing the polynomial
$\mathbf{g}= \lambda_1^5$.
Since $\lambda_1\neq 0$ we conclude the non-existence of  left-invariant  pp-wave metrics in the non-unimodular case.

	\begin{remark}\rm
Non-unimodular left-invariant metrics \eqref{eq:E-degenerate-+2} are Cotton-flat if and only if they are locally symmetric.
		Indeed, considering the same polynomial ring as above and computing a Gröbner basis for the ideal
		$\mathcal{I}= \langle \lambda_1^2 \Cotton_{ijk} \cup \{ (\gamma_1+\gamma_2)\varphi-1 \} \rangle$ we obtain
		46 polynomials, among which we have
		$\mathbf{g}_1=\lambda_1^3\lambda_2$, 
		$\mathbf{g}_2=(\gamma_1-\gamma_2)\lambda_1^4$,
		$\mathbf{g}_3=\gamma_3 \lambda_1^3$
		and
		$\mathbf{g}_4= \gamma_4 \lambda_1^3$. 
		Since $\lambda_1\neq 0$ it follows that  $\lambda_2=\gamma_3=\gamma_4=0$ and $\gamma_1=\gamma_2\neq 0$, and a direct calculation shows that the   space is locally symmetric.  
	\end{remark}

\subsubsection{\bf The induced metric on the derived algebra $\bm{\mathfrak{g}'_3}$ is degenerate}\label{EE-se:degenerate}

There exists a pseudo-orthonormal basis $\{u_1,u_2,u_3,u_4\}$ of the Lie algebra $\mathfrak{g}=\mathfrak{g}_3\rtimes\mathbb{R}$ with $\langle u_1,u_1\rangle=\langle u_2,u_2\rangle=\langle u_3,u_4\rangle=1$ so that $\mathfrak{g}_3=\operatorname{span}\{ u_1,u_2,u_3\}$ and the derived algebra $\mathfrak{g}'_3=\operatorname{span}\{ u_1,u_3\}$.
In this situation, using that $\mathfrak{g}_3$ is unimodular and $\mathfrak{g}'_3$ is Abelian, the Jacobi identity leads to two-different situations which we consider as follows.

\subsubsection*{\underline{Case 1}}
There exists a pseudo-orthonormal basis $\{u_1,u_2,u_3,u_4\}$ of the Lie algebra, with $\langle u_1,u_1\rangle=\langle u_2,u_2\rangle=\langle u_3,u_4\rangle=1$, so that
the left-invariant metrics are determined by
\begin{equation}\label{eq:E-degenerate-01}
\begin{array}{lll}
[u_1,u_2]\!=\!\lambda_1 u_3, &
[u_2,u_3]\!=\!\lambda_2 u_1,
\\ \noalign{\medskip}
[u_1,u_4]\!=\!\gamma_1 u_1 \!+\! \gamma_2 u_3, & 
[u_2,u_4]\!=\!\gamma_3 u_1\!+\!\gamma_4 u_3, & 
[u_3,u_4]\!=\!-\frac{\gamma_2\lambda_2}{\lambda_1} u_1\! +\! \gamma_1 u_3,
\end{array}
\end{equation}
where $\lambda_1\lambda_2\neq 0$ and  $\gamma_i\in\mathbb{R}$.

We work in the polynomial ring 
$\mathbb{R}[\gamma_1,\gamma_2,\gamma_3,\gamma_4,\lambda_1,\lambda_2]$.
Computing a Gröbner basis for the ideal $\mathcal{I}$ generated by 
$\lambda_1^4\RicciRiccicomp{i}{j} \cup \{\lambda_1^2 \tau\}$ (where we have multiplied by the powers of $\lambda_1$   to remove the variables in denominators) we get $24$ polynomials, one of them being
$\mathbf{g} = \lambda_1^5\lambda_2^3$.
Since $\lambda_1\lambda_2\neq 0$, we conclude that left-invariant metrics in this case are never pp-waves.

	\begin{remark}\rm
Left-invariant metrics \eqref{eq:E-degenerate-01} are never Cotton-flat. Indeed, considering the components of the Cotton tensor
		\[
		\Cotton_{132} = \tfrac{(\lambda_1^3-4\gamma_3\lambda_1^2+2\gamma_2^2\lambda_2)\lambda_2^2}{4\lambda_1^2},\quad
		\Cotton_{231} = \tfrac{(\lambda_1^3-8\gamma_3\lambda_1^2+4\gamma_2^2\lambda_2)\lambda_2^2}{4\lambda_1^2},
		\] 
		one has
		$2\Cotton_{132}-\Cotton_{231} = \frac{1}{4}\lambda_1\lambda_2^2$, which never vanishes since
		$\lambda_1\lambda_2\neq 0$.
	\end{remark}

\subsubsection*{\underline{Case 2}}
There exists a pseudo-orthonormal basis $\{u_1,u_2,u_3,u_4\}$ of the Lie algebra, with $\langle u_1,u_1\rangle=\langle u_2,u_2\rangle=\langle u_3,u_4\rangle=1$, so that 
\begin{equation}\label{eq:E-degenerate-02}
\begin{array}{l}
[u_1,u_2]= \lambda_1 u_1+\lambda_2 u_3,  \quad
[u_2,u_3]= \lambda_3 u_1+\lambda_1 u_3, \quad
[u_2,u_4]= \gamma_2 u_1 + \gamma_3 u_3,

\\ \noalign{\medskip}

[u_1,u_4]= \gamma_1 u_1 + \frac{(\gamma_1-\gamma_4)\lambda_2}{2\lambda_1} u_3, \quad
[u_3,u_4]= - \frac{(\gamma_1-\gamma_4)\lambda_3}{2\lambda_1} u_1+\gamma_4 u_3,  
\end{array}
\end{equation}
where $\lambda_1\neq 0$, $\lambda_1^2-\lambda_2\lambda_3\neq 0$ and  $\gamma_i\in\mathbb{R}$.

We take an auxiliary variable $\lambda'_1$ to indicate that $\lambda_1\neq 0$ by means of $\lambda_1\lambda'_1-1$.
We consider  the ideal $\mathcal{I}=\langle \lambda_1^4 \RicciRiccicomp{i}{j} \cup \{\lambda_1^2 \tau,  \lambda_1\lambda'_1-1\}\rangle$ in the polynomial ring
$\mathbb{R}[\gamma_1,\gamma_2,\gamma_3,\gamma_4,\lambda'_1,\lambda_1,\lambda_2,\lambda_3]$  (where as in the previous cases we have multiplied by the powers of $\lambda_1$ to remove the variables in denominators) and compute a Gröbner basis for this ideal. Thus we get $8$ polynomials, among which we find
$\mathbf{g} = \lambda_1^2-\lambda_2\lambda_3$.
Since $\lambda_1^2-\lambda_2\lambda_3\neq 0$, we conclude that left-invariant metrics in this case are never pp-waves.

 	\begin{remark}\rm
 Left-invariant metrics \eqref{eq:E-degenerate-02} are never Cotton-flat.
 		To show it, we consider the ideal 
 		$\mathcal{I}= \langle \lambda_1^3 \Cotton_{ijk} \cup \{ \lambda_1\lambda'_1-1 \} \rangle$ 
 		in the same polynomial ring as above and compute a Gröbner basis, obtaining 6 polynomials among which one has
 		$\mathbf{g}=\lambda_1^2-\lambda_2\lambda_3$, which does not vanish. 
 	\end{remark}

\section{Direct extensions of the non-solvable Lie groups}\label{se:SLSU}

We proceed as in the previous sections considering all the possible left-invariant Lorentzian metrics on the non-solvable four-dimensional Lie groups $\widetilde{SL}(2,\mathbb{R})\times\mathbb{R}$ and $SU(2)\times\mathbb{R}$. The different cases corresponding to a Riemannian, Lorentzian or degenerate metric on $\widetilde{SL}(2,\mathbb{R})$ or $SU(2)$ are considered separately in Sections~\ref{se:91}, \ref{se:92} and  \ref{se:93}, respectively. 
As a consequence of the analysis below we show that left-invariant metrics in this setting are neither pp-waves nor non-trivial Cotton-flat. Even more, Gröbner bases are again a fundamental tool to obtain the following result:

\begin{theorem}\label{SL-SU-th: pp-waves}
There is no left-invariant Lorentz metric on $\widetilde{SL}(2,\mathbb{R})\times\mathbb{R}$ or $SU(2)\times\mathbb{R}$ with two-step nilpotent Ricci operator. 

Moreover, a left-invariant Lorentz metric on $\widetilde{SL}(2,\mathbb{R})\times\mathbb{R}$ or $SU(2)\times\mathbb{R}$ has harmonic curvature if and only if it is locally conformally flat or a product metric on $\widetilde{SL}(2,\mathbb{R})\times\mathbb{R}$ determined by a locally conformally flat metric on $\widetilde{SL}(2,\mathbb{R})$.
\end{theorem}

\subsection{Direct extensions with Riemannian Lie groups $\bm{\widetilde{SL}(2,} \pmb{\mathbb{R}}\bm{)}$ or $\bm{SU(2)}$} \label{se:91}
Let $G=G_3\times \mathbb{R}$ with $G_3=\widetilde{SL}(2, \mathbb{R})$ or $G_3=SU(2)$ so that the restriction of the metric to $G_3$ is positive definite. Hence, there exists  an orthonormal basis $\{e_1,e_2,e_3,e_4\}$ of the Lie algebra, with $e_4$ timelike,  so that    
\begin{equation}\label{S-Riemann}
\begin{array}{l}
[e_1,e_2]=\lambda_3 e_3,\qquad
[e_1,e_3]=-\lambda_2 e_2, \qquad
[e_2,e_3]=\lambda_1 e_1,

\\
\noalign{\medskip}

[e_1,e_4]=\gamma_1\lambda_2 e_2+\gamma_2\lambda_3 e_3, \qquad 
[e_2,e_4]=-\gamma_1\lambda_1 e_1+\gamma_3\lambda_3 e_3, 

\\
\noalign{\medskip}

[e_3,e_4]=-\gamma_2\lambda_1 e_1-\gamma_3\lambda_2 e_2,
\end{array}
\end{equation}
where $\lambda_1\lambda_2\lambda_3\neq 0$ and  $\gamma_i\in\mathbb{R}$. 

Let $\mathcal{I}=\langle \RicciRiccicomp{i}{j} \cup \{\tau\} \rangle
\subset \mathbb{R}[\gamma_1,\gamma_2,\gamma_3,\lambda_1,\lambda_2,\lambda_3]$. Computing a Gröbner basis for this ideal  we get $212$ polynomials, one of them being
$\mathbf{g} = \lambda_1^3\lambda_2^3\lambda_3^3$.
Since $\lambda_1\lambda_2\lambda_3\neq 0$  we conclude that $\Ricci^2\neq 0$ and therefore
there is no left-invariant pp-wave metric in this case.

	\begin{remark}\label{re: SL-SU Cotton-flat}\rm
		To analyze the vanishing of the Cotton tensor for a left-invariant metric \eqref{S-Riemann}, we use an auxiliary variable $\varphi$ to indicate that $\lambda_1$, $\lambda_2$ and $\lambda_3$ are non-zero by means of the polynomial $\lambda_1\lambda_2\lambda_3 \varphi - 1$. Now,  we consider the ideal $\mathcal{I}=\langle \Cotton_{ijk}\cup \{ \lambda_1\lambda_2\lambda_3 \varphi - 1\} \rangle$ in the polynomial ring $\mathbb{R}[\varphi,\gamma_1,\gamma_2,\gamma_3,\lambda_1,\lambda_2,\lambda_3]$
		and compute a Gröbner basis, obtaining $12$ polynomials, two of them being $\mathbf{g}_1=(\lambda_2-\lambda_3)(\lambda_2^2+\lambda_3^2+\lambda_2\lambda_3)$ and $\mathbf{g}_2= \lambda_1\lambda_2-\lambda_3^2$. Hence, $\lambda_1=\lambda_2=\lambda_3$, from where it follows that the metric is locally conformally flat locally isometric to $\mathbb{S}^3\times\mathbb{R}$.
	\end{remark}

\subsection{Direct extensions with Lorentzian Lie groups $\bm{\widetilde{SL}(2,} \pmb{\mathbb{R}}\bm{)}$ or $\bm{SU(2)}$} \label{se:92}
We consider the different possibilities for the Jordan normal form of the structure operator.

\subsubsection{\bf Diagonalizable structure operator}\label{sl-Ia}
There exists an orthonormal basis $\{e_1,e_2,e_3,e_4\}$ of the Lie algebra, with $e_3$ timelike,  so that  
\begin{equation}\label{S-Lorentz-Ia}
\begin{array}{l}
[e_1,e_2]=-\lambda_3 e_3,\qquad
[e_1,e_3]=-\lambda_2 e_2,  \qquad
[e_2,e_3]=\lambda_1 e_1,

\\
\noalign{\medskip}
[e_1,e_4]=\gamma_1\lambda_2 e_2+\gamma_2\lambda_3 e_3,\qquad
[e_2,e_4]=-\gamma_1\lambda_1 e_1+\gamma_3\lambda_3 e_3, 

\\
\noalign{\medskip}

[e_3,e_4]=\gamma_2\lambda_1 e_1+\gamma_3\lambda_2 e_2,
\end{array}
\end{equation}
where $\lambda_1\lambda_2\lambda_3\neq 0$  and  $\gamma_i\in\mathbb{R}$. 

We proceed exactly as in the previous case and conclude that $\Ricci^2\neq 0$. Therefore
there is no left-invariant pp-wave metric in this case.

	\begin{remark}\rm
		We proceed exactly as in Remark~\ref{re: SL-SU Cotton-flat} to get that any Cotton-flat left-invariant metric \eqref{S-Lorentz-Ia} must be locally conformally flat and locally isometric to $\mathbb{S}^3_1\times\mathbb{R}$.
	\end{remark}

\subsubsection{\bf Structure operator $\bm{L}$ with a complex eigenvalue}\label{sl-complex}
In this case  there exists an orthonormal basis $\{e_1,e_2,e_3,e_4\}$ of the Lie algebra, with $e_3$ timelike,  such that the Lie algebra structure is given by
\[
\begin{array}{l}
[e_1,e_2]=-\beta e_2-\alpha e_3,\quad\qquad\qquad\,\,\,\,
[e_1,e_3]=-\alpha e_2 + \beta e_3,\quad
[e_2,e_3]=\lambda e_1, 
\\
\noalign{\medskip}
[e_1,e_4]=(\alpha^2+\beta^2)(\gamma_1 e_2+\gamma_2 e_3), 
\quad
[e_2,e_4]=-(\gamma_1\alpha-\gamma_2\beta)\lambda e_1 +  \gamma_3 \beta e_2 +  \gamma_3 \alpha e_3, 
\\
\noalign{\medskip}
[e_3,e_4]=(\gamma_2\alpha+\gamma_1\beta)\lambda e_1 + \gamma_3 \alpha e_2 -   \gamma_3 \beta e_3 ,
\end{array}
\]
where    $\beta\lambda\neq 0$ and $\alpha, \gamma_i\in\mathbb{R}$.
In this case the three-dimensional unimodular Lie algebra corresponds to $\mathfrak{sl}(2,\mathbb{R})$.

In the polynomial ring $\mathbb{R}[\lambda,\alpha,\beta,\gamma_1,\gamma_2,\gamma_3]$ we consider the ideal $\mathcal{I}$ generated by $\RicciRiccicomp{i}{j} \cup \{\tau\}$. Computing a Gröbner basis for this ideal we obtain a set of $108$ elements containing the polynomial
$\mathbf{g} =  (\gamma_3^2+1)^3 (\alpha^2+\beta^2) \beta^3$,
which does not vanish since $\beta\neq 0$. Hence we conclude that $\Ricci^2\neq 0$ and therefore left-invariant metrics in this case are never pp-waves.

	\begin{remark}\rm
		Left-invariant Cotton-flat metrics on $\widetilde{SL}(2,\mathbb{R})\times\mathbb{R}$  whose structure operator has a complex eigenvalue reduce to the product metric, where the metric on $\widetilde{SL}(2,\mathbb{R})$ is the locally conformally flat one described in Section~\ref{section-3D-Cotton flat}.
	Indeed, considering the same polynomial ring as above and computing a Gröbner basis for the ideal $\mathcal{I}=\langle \Cotton_{ijk}\rangle$, one obtains $44$ polynomials, among which we find
		$\mathbf{g}_1=(\gamma_1^2+\gamma_2^2)(\gamma_3^2+1)\beta^5$ and
		$\mathbf{g}_2=\gamma_3(\gamma_3^2+1)^2\beta^5$. Hence, $\gamma_1=\gamma_2=\gamma_3$ and the left-invariant metric is the product one, from where it follows that $\widetilde{SL}(2,\mathbb{R})\times\mathbb{R}$ is Cotton-flat if and only if so is $\widetilde{SL}(2,\mathbb{R})$ with the induced Lorentzian metric, which corresponds to the metric (1) in Section~\ref{section-3D-Cotton flat}. 
	\end{remark}

\subsubsection{\bf Structure operator of type II}\label{sl-II}
In this case the left-invariant metrics are described by
\begin{equation}\label{S-lorentz-II}
\begin{array}{l}
[u_1,u_2]= \lambda_2 u_3 ,\qquad
[u_1,u_3]=-\lambda_1 u_1-\ve u_2,\qquad
[u_2,u_3]=\lambda_1 u_2,

\\
\noalign{\medskip}
[u_1,u_4]=\gamma_1\lambda_1 u_1 + \ve\gamma_1 u_2 + \gamma_2\lambda_2 u_3, 
\quad
[u_2,u_4]=-\gamma_1\lambda_1 u_2 + \gamma_3 \lambda_2 u_3, 
\\
\noalign{\medskip}
[u_3,u_4]=-\gamma_3\lambda_1 u_1 - (\gamma_2\lambda_1+\ve\gamma_3)u_2 ,
\end{array}
\end{equation}
where   $\varepsilon^2=1$, $\lambda_1\lambda_2\neq 0$,  $\gamma_i\in\mathbb{R}$ and
$\{u_1,u_2,u_3,u_4\}$  is a pseudo-orthonormal basis of the Lie algebra with $\langle u_1,u_2\rangle=\langle u_3,u_3\rangle=\langle u_4,u_4\rangle=1$.

Computing a Gröbner basis for the ideal $\mathcal{I} =\langle \RicciRiccicomp{i}{j} \cup \{\tau,\ve^2-1\}\rangle$ in the polynomial ring  
$\mathbb{R}[\ve,\gamma_1,\gamma_2,\gamma_3,\lambda_1,\lambda_2]$  we get $41$ polynomials, one of them being
$\mathbf{g}=\lambda_2^5$.
Since $\mathbf{g}\neq 0$, it follows that  $\Ricci^2\neq 0$ and therefore there is no left-invariant pp-wave metric in this setting.

	\begin{remark}\rm
Left-invariant metrics \eqref{S-lorentz-II} are never Cotton-flat. Indeed, working in the same polynomial ring as above, we  compute a Gröbner basis for the ideal $\mathcal{I}=\langle \Cotton_{ijk}\cup \{\ve^2-1\} \rangle$ to get  $27$ polynomials,
		one of them being 
		$\mathbf{g}=\lambda_2^4$, which does not vanish.
	\end{remark}

\subsubsection{\bf Structure operator of type III}\label{sl-III}
There exists a pseudo-orthonormal basis $\{u_1,u_2,u_3,u_4\}$ of the Lie algebra, with $\langle u_1,u_2\rangle=\langle u_3,u_3\rangle=\langle u_4,u_4\rangle=1$,  such that the corresponding Lie brackets are determined by 
\[
\begin{array}{l}
[u_1,u_2]= u_1 + \lambda u_3 ,\qquad\qquad
[u_1,u_3]=-\lambda u_1 ,\qquad
[u_2,u_3]=\lambda u_2+u_3,

\\
\noalign{\medskip}
[u_1,u_4]=\gamma_1\lambda u_1 + \gamma_2 \lambda^2 u_3,
%
\quad
[u_2,u_4]=\gamma_3 u_1 -(\gamma_1-\gamma_2)\lambda u_2 -(\gamma_1-\gamma_2-\gamma_3\lambda)u_3 ,

\\
\noalign{\medskip}

[u_3,u_4]=-\gamma_3\lambda u_1-\gamma_2\lambda^2 u_2-\gamma_2\lambda u_3 ,
\end{array}
\]
where $\lambda\neq 0$ and   $\gamma_i\in\mathbb{R}$.  
A direct calculation shows that
\[
\Riccicomp{1}{1} = \Riccicomp{2}{2} =
		\tfrac{-1}{2}(\gamma_2^2-\gamma_1\gamma_2+1)\lambda^2,
\,\,\,
\Riccicomp{3}{3} = \tfrac{1}{2}(2\gamma_2^2-2\gamma_1\gamma_2-1)\lambda^2,
\,\,\,
\Riccicomp{4}{4} = \tfrac{-3}{2}\gamma_2^2\lambda^2,
\]
so that the scalar curvature $\tau=-\frac{3}{2}(\gamma_2^2+1)\lambda^2$ does not vanish. Hence the Ricci operator cannot be nilpotent and left-invariant metrics in this case are never pp-waves.

	\begin{remark}\rm
Considering the components of the Cotton tensor
		\[
		\Cotton_{121} = -\tfrac{9}{4}\gamma_2^2\lambda^4,\quad
		\Cotton_{122} = \tfrac{3}{4}\left( \gamma_1^2+2 -\gamma_2 (2\gamma_2+5\gamma_3\lambda) \right)\lambda^2,
		\]
		one has that $\gamma_2=0$ and thus $\Cotton_{122}=\frac{3}{4}(\gamma_1^2+2)\lambda^2$, which does not vanish. Hence there is no Cotton-flat left-invariant metrics  with structure operator of type III.
	\end{remark}

\subsection{Direct extensions with degenerate Lie groups $\bm{\widetilde{SL}(2,} \pmb{\mathbb{R}}\bm{)}$ or $\bm{SU(2)}$} \label{se:93}

Next we assume that the restriction of the left-invariant metric to $\widetilde{SL}(2,\mathbb{R})$ or $SU(2)$ is degenerate.
Let $u$ be such that $\span\{ u\}$ is degenerate. We consider separately the cases of $\ad_{u}$ having two purely imaginary complex eigenvalues,   two non-zero real opposite eigenvalues, or being   three-step nilpotent.

\subsubsection{\bf  $\bm{{\ad_{u}}}$ has complex eigenvalues}\label{se:complex}

There exists a basis  $\{ v_1,v_2,v_3,v_4\}$ of the Lie algebra, with $\langle v_1,v_1\rangle=\langle v_2,v_2\rangle=\langle v_3,v_4\rangle=1$, $\langle v_1,v_2\rangle=\kappa$ with $\kappa^2<1$,   such that the Lie brackets are  given by (up to an isomorphic homothety)
\begin{equation}\label{S-degenerate-complex}
\begin{array}{lll}
\![v_1,v_2]=v_3, &
\!\![v_1,v_3]=\beta\lambda^2 v_2, &
\!\![v_1,v_4]= \gamma_1\lambda^2 v_2+\gamma_2 v_3,	
\\
\noalign{\medskip}
\![v_2,v_3]=-\beta v_1, &
\!\![v_2,v_4]= -\gamma_1 v_1 + \gamma_3 v_3, &
\!\![v_3,v_4]= \gamma_2\beta v_1+\gamma_3\beta\lambda^2 v_2 ,
\end{array}
\end{equation}
where $\beta\lambda\neq 0$ and  $\gamma_i\in\mathbb{R}$.

We introduce auxiliary variables $\lambda'$, $\beta'$ and $\varphi$ to express that $\lambda$, $\beta$ and $\kappa^2-1$ are non-zero by means of the polynomials $\lambda\lambda'-1$, $\beta\beta'-1$ and $(\kappa^2-1)\varphi-1$.
Let $\mathcal{I}=\langle  (\kappa^2-1)^2 \RicciRiccicomp{i}{j} \cup \{(\kappa^2-1)\beta^{-1} \tau, 
\lambda\lambda'-1,
\beta\beta'-1,
(\kappa^2-1)\varphi-1\}\rangle$
in the polynomial ring 
$\mathbb{R}[\lambda,\kappa,\gamma_1,\gamma_2,\gamma_3,\beta,\lambda',\beta',\varphi]$,
where we have multiplied by the powers of $\kappa^2-1$   to remove the variables in denominators.
Computing a Gröbner basis for this ideal  it reduces to $\{ 1\}$. Therefore  $\Ricci^2\neq 0$, which shows that    this case does not support
any left-invariant pp-wave metric.

\begin{remark}\rm
	Left-invariant metrics \eqref{S-degenerate-complex} are never Cotton-flat. Indeed, it follows from the components
	$\Cotton_{133}=-\frac{1}{4}\beta^3\lambda^2(\lambda^2+1)\gamma_2$ and $\Cotton_{233}=-\frac{1}{4}\beta^3\lambda^2(\lambda^2+1)\gamma_3$ that $\gamma_2=\gamma_3=0$. Then, one has the component 
	$
	\Cotton_{124}= \frac{\beta(\lambda^2+1)
		\left( 
		\left(
		8\kappa^2\lambda^2+2(\lambda^2-1)^2
		\right)\gamma_1^2+1\right)}{2(\kappa^2-1)}$, 
	which does not vanish.
\end{remark}

\black

\subsubsection{\bf $\bm{{\ad_{u}}}$ has non-zero real eigenvalues}\label{se:real}

There exists a basis $\{ v_1,v_2,v_3,v_4\}$ of the Lie algebra, with $\langle v_1,v_1\rangle=\langle v_2,v_2\rangle=\langle v_3,v_4\rangle=1$, $\langle v_1,v_2\rangle=\kappa$ with $\kappa^2<1$, so that   
\begin{equation}\label{S-degenerate-real}
\begin{array}{lll}
	
	\![v_1,v_2]=v_3, &
	\![v_1,v_3]=\lambda v_1, &
	\![v_1,v_4]= \gamma_1 v_1+\gamma_2 v_3,	
	
	\\
	\noalign{\medskip}
	
	\![v_2,v_3]=-\lambda v_2, &
	\![v_2,v_4]= -\gamma_1 v_2 + \gamma_3 v_3, &
	\![v_3,v_4]= \gamma_3\lambda v_1+\gamma_2\lambda v_2,
\end{array}
\end{equation}
where $\lambda\neq 0$ and  $\gamma_i\in\mathbb{R}$. 

We proceed as in the previous case, multiplying  by the powers of $\kappa^2-1$   to remove the variables in denominators, when we consider $\mathcal{I}=\langle 
(\kappa^2-1)^2 \RicciRiccicomp{i}{j} \cup \{(\kappa^2-1) \tau\} \rangle$ 
in the polynomial ring
$\mathbb{R}[\gamma_1,\gamma_2,\gamma_3,\kappa, \lambda]$.
Computing a Gröbner basis for this ideal  we get $3$ polynomials, one of them being
$\mathbf{g} = \lambda^2$.
Since $\lambda\neq 0$ we have that $\Ricci^2\neq 0$, which  shows the   non-existence of left-invariant pp-wave metrics in this case.

	\begin{remark}\rm
	Left-invariant metrics \eqref{S-degenerate-real} have harmonic curvature if and only it they are locally conformally flat. Indeed, considering the same polynomial ring as above and computing a Gröbner basis for the ideal 
		$\mathcal{I}=\langle 
		(\kappa^2-1) \Cotton_{ijk}  \rangle$  one gets 11 polynomials, among which we find
		$\mathbf{g}_1= \gamma_1 \lambda^2$,
		$\mathbf{g}_2= \gamma_2 \lambda^2$,
		$\mathbf{g}_3= \gamma_3 \lambda^2$ and
		$\mathbf{g}_4= \kappa \lambda^2$. Since $\lambda\neq 0$, it follows that  $\gamma_1=\gamma_2=\gamma_3=\kappa=0$ and a direct calculation shows that the left-invariant metric is locally conformally flat but not locally symmetric. 
	\end{remark}

\subsubsection{\bf $\bm{{\ad_{u}}}$ is three-step nilpotent}\label{se:nilpotent}
The left-invariant metrics   are determined, with respect to a basis  $\{v_1,v_2,v_3,v_4\}$ of the Lie algebra  with $\langle v_1,v_1\rangle=\langle v_2,v_2\rangle=\langle v_3,v_4\rangle=1$ and $\langle v_1,v_2\rangle=\kappa$ with $\kappa^2<1$, by the Lie brackets
\[
\begin{array}{lll}
	\![v_1,v_2]=\alpha v_1+\beta v_3, &
	\!\![v_1,v_3]=-v_2, &
	\!\! [v_1,v_4]= \gamma_1 v_1+\gamma_2 v_2+\frac{\beta\gamma_1}{\alpha}v_3,	
	\\
	\noalign{\medskip}
	\![v_2,v_3]=\alpha v_3, &
	\!\! [v_3,v_4]= \gamma_3 v_2-\gamma_1 v_3 , &
	\!\! [v_2,v_4]= -\alpha\gamma_3 v_1 - (\alpha\gamma_2+\beta\gamma_3) v_3,
\end{array}
\]
where $\alpha\neq 0$, $\kappa^2\neq 1$ and $\beta,\gamma_i\in\mathbb{R}$.

We consider  
$\mathcal{I}=\langle 
(\kappa^2-1)^2 \alpha^2 \RicciRiccicomp{i}{j} \cup \{(\kappa^2-1) \tau\} \rangle
\subset
\mathbb{R}[\beta,\gamma_1,\gamma_2,\gamma_3,\alpha,\kappa]$, where once again  we have multiplied  by the powers of $\alpha$ and $\kappa^2-1$   to remove the variables in denominators. Computing a Gröbner basis for this ideal we obtain 10 polynomials, among which we find
$
\mathbf{g}= (\kappa^2-1)\alpha^5$.
Since $\alpha\neq 0$ and $\kappa^2\neq 1$, it follows that $\Ricci^2\neq 0$ and we conclude that there is no left-invariant pp-wave metric in this case.

	\begin{remark}\rm
		A direct calculation shows that $\Cotton_{233}=\frac{\alpha}{4(1-\kappa^2)}\neq 0$, which shows that there is no Cotton-flat left-invariant metric in this case.
	\end{remark}

\end{document}